\documentclass[12pt,a4paper]{article}
\pdfoutput=1

\usepackage[T1]{fontenc}
\usepackage[UKenglish]{isodate}
\cleanlookdateon

\usepackage[style=numeric,maxcitenames=99, maxbibnames=99, giveninits=true]{biblatex} 
\usepackage{amsmath,amssymb,mathtools,dsfont,microtype} 
\usepackage[dvipsnames]{xcolor} 
\usepackage{pdfcolfoot} 
\usepackage[noBBpl]{mathpazo} 
\usepackage{bbm}
\usepackage[labelsep=period,labelfont=bf,justification=centering]{caption}
\usepackage{float,graphicx,subcaption}
\usepackage{booktabs,makecell}
\usepackage{enumitem,mdwlist}
\setlist[itemize]{topsep=0ex,itemsep=0ex,parsep=0.4ex}
\setlist[enumerate]{topsep=0ex,itemsep=0ex,parsep=0.4ex}

\usepackage{parskip,fullpage}
\usepackage{standalone}

\usepackage[linktoc=all,hidelinks,colorlinks,unicode=true]{hyperref} 
\usepackage[capitalise,compress,nameinlink,noabbrev]{cleveref} 

\hypersetup{linkcolor={blue!70!black}, citecolor={green!70!black}, urlcolor={blue!70!black}}

\usepackage{amsthm}

\theoremstyle{plain}
\newtheorem{theorem}{Theorem}[section]

\newtheorem{conjecture}[theorem]{Conjecture}
\newtheorem{lemma}[theorem]{Lemma}

\newtheorem{proposition}[theorem]{Proposition}

\theoremstyle{definition}

\crefname{definition}{Definition}{Definitions}
\crefname{theorem}{Theorem}{Theorems}
\crefname{lemma}{Lemma}{Lemmas}
\crefname{conjecture}{Conjecture}{Conjectures}
\crefname{claim}{Claim}{Claims}
\crefname{corollary}{Corollary}{Corollaries}
\crefname{observation}{Observation}{Observations}
\crefname{proposition}{Proposition}{Propositions}
\crefname{question}{Question}{Questions}
\crefname{remark}{Remark}{Remarks}
\crefname{example}{Example}{Examples}
\crefname{section}{Section}{Sections}
\crefname{subsection}{\S}{\S\S}

\newcommand*{\cA}{\mathcal{A}}

\newcommand*{\cS}{\mathcal{S}}

\renewcommand{\epsilon}{\varepsilon}
\renewcommand{\ge}{\geqslant}
\renewcommand{\le}{\leqslant}
\renewcommand{\geq}{\geqslant}
\renewcommand{\leq}{\leqslant}

\renewcommand{\emptyset}{\varnothing}
\renewcommand{\subset}{\subseteq}

\newcommand{\defn}[1]{{\color{Maroon}\ifmmode#1\else\emph{#1}\fi}}

\title{Pairwise incomparable and uncomplemented families }

\date{\today}
\author{Julia B\"ottcher\footnotemark[1] \qquad  Pedram Noohi\footnotemark[2] \\ Jozef Skokan\footnotemark[3] \qquad John Talbot\footnotemark[4] 
}
\begin{document}

\maketitle
\newcommand{\pedram}[1]{{\color{blue} P: #1}}

\begin{abstract}
In 1976 Hilton proved that whenever families $\mathcal{A}_1, \mathcal{A}_2\subseteq2^{[n]}$ are uncomplemented (i.e. $\bar{A}\not\in \mathcal{A}_i$ for $A\in \mathcal{A}_i$) and pairwise incomparable (i.e. if $A\in \mathcal{A}_i$ and $B\in \mathcal{A}_j$ for $i\neq j$, then $A\not\subseteq B$ and $B\not\subseteq A$), then necessarily  $|\mathcal{A}_1|+|\mathcal{A}_2|\leq 2^{n-1}$. He then conjectured that extending this to $k$ families $\mathcal{A}_1\ldots ,\mathcal{A}_k$ would lead to $\sum_{i=1}^{k}|\mathcal{A}_i|\leq 2^{n-1}$. We prove this conjecture, and show that if at least three of the families are non-empty then in fact $\sum_{i=1}^{k}|\mathcal{A}_i|\leq 2^{n-1}-1$, which is tight whenever $3\leq k\leq n$. Finally we show, again for $k \ge 3$, that if for every $i\not=j$ there is a complementary pair between $\mathcal{A}_i$ and $\mathcal{A}_j$, then $\sum_{i=1}^k|\mathcal{A}_i|\leq 2^{n-1}-2$. 
\end{abstract}
\renewcommand{\thefootnote}{\fnsymbol{footnote}} 

\footnotetext[1]{Department of Mathematics, London School of Economics, UK (\textsf{\href{mailto:j.Boettcher@lse.ac.uk}{j.Boettcher@lse.ac.uk}})}

\footnotetext[2]{Department of Mathematics, London School of Economics, UK (\textsf{\href{mailto:p.noohi@lse.ac.uk}{p.noohi@lse.ac.uk}})}

\footnotetext[4]{Department of Mathematics, University College London, UK (\textsf{\href{mailto:j.talbot@ucl.ac.uk}{j.talbot@ucl.ac.uk}})}

\footnotetext[3]{Department of Mathematics, London School of Economics, UK (\textsf{\href{mailto:j.skokan@lse.ac.uk}{j.skokan@lse.ac.uk}})}

\renewcommand{\thefootnote}{\arabic{footnote}} 
\section{Introduction}
 Families $\mathcal{A}_1,\mathcal{A}_2,\dots,\mathcal{A}_k\subseteq 2^{[n]}$ are called \emph{pairwise incomparable} if for $i,j\leq k$ and $i\not=j$  we have that for $A\in \mathcal{A}_i$ and $B\in\mathcal{A}_j$ that $A\not\subseteq B$ and $B\not\subseteq A$.
A family $\mathcal{A}\subseteq 2^{[n]}$ is an \emph{anti-chain} if for any $A_1,A_2\in\mathcal{A}$ it is the case that $A_1\not\subseteq A_2$ and $A_2\not \subseteq A_1$, and a classical result in Extremal Set Theory of Sperner \cite{sperner} says that an \text{anti-chain} has size at most $\binom{n}{\lceil n/2 \rceil}$. 
An extension of the notion of an \text{anti-chain} would be taking $k,n\in \mathbb{N}$ and $\mathcal{A}_1,\mathcal{A}_2,...,\mathcal{A}_k\subseteq 2^{[n]}$ such that the families are pairwise incomparable which is at times called a cross-Sperner system. 

Work on bounding the sizes of sums of pairwise incomparable families has been studied before, most notably by Seymour \cite{Seymour1973} resolving a Conjecture of Hilton, showing that if $\mathcal{A}\subseteq2^{[n]}$ and $\mathcal{B}\subseteq 2^{[n]}$ are pairwise incomparable then $|\mathcal{A}|^{1/2}+|\mathcal{B}|^{1/2}\leq 2^{n/2}$. Then in $2011$ Gerbner, Lemons, Palmer, Patkós and Szécsi \cite{GerbnerEtAl2012} showed that for sufficiently large $n$ and $\mathcal{A},\mathcal{B}\subseteq2^{[n]}$ pairwise incomparable and non-empty, that $|\mathcal{A}|+|\mathcal{B}|\leq 2^n-2^{\lceil n/2\rceil}-2^{\lfloor n/2\rfloor}+2$ which gives the same upper bound as Seymour's for even $n$, but is stronger for $n$ odd. Letting $$\sigma(n,k)=\max\left\{\sum_{i=1}^k|\mathcal{A}_i|:\mathcal{A}_1,...,\mathcal{A}_k\subseteq 2^{[n]} \text{ are pairwise incomparable and non-empty}\right\}$$ Behague, Kuperus, Morrison and Wright \cite{BehagueEtAl2024} showed that for $n\geq 2k$ that
\begin{equation*}
2^n-\frac{3}{\sqrt{2}}\sqrt{2^nk}+2(k-1)\leq \sigma(n,k)\leq 2^n-2\sqrt{2^n(k-1)}+2(k-1).
\end{equation*}
An example of another property studied for pairwise incomparable families is the following: For two pairwise incomparable families $\mathcal{A},\mathcal{B}\subseteq2^{[n]}$ define $I(\mathcal{A},\mathcal{B})=\{A\cap B:A\in\mathcal{A}\text{ and } B\in\mathcal{B}\}.$ Frankl and Wang \cite{FranklWang2023Intersections} initiated the study of $m(n)=\max\{|I(\mathcal{A},\mathcal{B})|:\mathcal{A}\subseteq2^{[n]}\text{ and }\mathcal{B}\subseteq2^{[n]}\text{ are pairwise incomparable}\}$ and Pan \cite{Pan2023CrossSperner} showed that $m(n)=2^n-2^{\lfloor n/2\rfloor}-2^{\lceil n/2\rceil}+1$.

A family $\mathcal{A}\subseteq 2^{[n]}$ is \emph{uncomplemented} if for any $A\in\mathcal{A}$ it is the case that $\overline{A}\not\in\mathcal{A}$, and so clearly an uncomplemented family has size at most $2^{n-1}$.
In  $1976$, Hilton \cite{Hilton1976} initiated the study of bounding the size of the sum of uncomplemented and pairwise incomparable families of sets, proving in a short paper that when $k=2$ that it is the case that $|\mathcal{A}_1|+|\mathcal{A}_2|\leq2^{n-1}$, which is tight by taking $x,y\in [n]$ and $x\not=y$ with $\mathcal{A}_1=\{A\cup\{x\}:A\subseteq [n]\setminus \{x,y\}\}$ and $\mathcal{A}_2=\{A\cup\{y\}:A\subseteq [n]\setminus \{x,y\}\}$. In a review of Hilton's paper \cite{KleitmanReview1977}, Kleitman reproved the result 
using the Harris-Kleitman inequalities (see Theorem \ref{Harris}) and, appearing in the same review, Hilton conjectured the following.
\begin{conjecture}
\label{conj}
For $k\in\mathbb{N}$, $k\geq 3$ and $\mathcal{A}_1,\dots,\mathcal{A}_k\subseteq2^{[n]}$ where each family is uncomplemented and the families are pairwise incomparable it is the case that,
\begin{equation*}\sum_{i=1}^k|\mathcal{A}_i|\leq 2^{n-1}.
\end{equation*}
\end{conjecture}
Partial progress on this conjecture was made first in $2001$ by Liu and Zhao \cite{LiuZhao2001} where they proved the conjectural bound under different assumptions, and the general bound of  $\displaystyle\sum_{i=1}^k|\mathcal{A}_i|\leq \tfrac{k2^n}{2\sqrt{k-1}+k}$ or $\displaystyle\sum_{i=1}^k|\mathcal{A}_i|<\tfrac{k2^{n-1}}{k-1}$. Then in $2005$, Peng and Zhao \cite{PengZhao2005} improved the bound to $\displaystyle\sum_{i=1}^k|\mathcal{A}_i|< (4-2\sqrt{2})2^{n-1}$ using a generalisation of Seymour's result to $k$ many families, and showing that the Conjecture holds under certain conditions on the sizes of a subset of the families. 

In this paper we first prove Conjecture \ref{conj} in Section \ref{firstsection}, and use the tools developed to prove Theorem \ref{Main} in Section \ref{secondsection}, which is a slight strengthening of the Conjecture. We do this by studying the connected components of the comparability graph of the family of sets.
\begin{theorem}
\label{Main}
 Let $\mathcal A_1,\dots\mathcal A_k\subseteq2^{[n]}$ be, uncomplemented and pairwise incomparable families, then if $k\geq 3$ and at least $3$ families are non-empty, it is the case that \begin{equation*}
 \sum_{i=1}^k|\mathcal{A}_i|\le2^{n-1}-1.
 \end{equation*}
\end{theorem}
We also show in Proposition \ref{tight} that Theorem \ref{Main} is tight when $3 \le k\leq n$.
\begin{proposition}
\label{tight}
For $3 \le k\leq n$ there exist non-empty families $\mathcal{A}_1,\dots,\mathcal{A}_k$ that are uncomplemented and pairwise incomparable such that \begin{equation*}
\sum^k_{i=1}|\mathcal{A}_i|=2^{n-1}-1.
\end{equation*}
\end{proposition}
In the last section of the paper we then impose the condition that between any $2$ families $\mathcal{A}_i$ and $\mathcal{A}_j$ for $i\not=j$ it is the case that there exists some $A\in \mathcal{A}_i$ such that $\overline A\in \mathcal{A}_j$. We call such collection of families \emph{non-degenerate} in the sense that no two families can be merged together. Given this definition, Section \ref{lastsection} is devoted to proving the following Theorem.

\begin{theorem}
\label{3rd}
 Let $k \ge 3$ and $\mathcal A_1,\dots\mathcal A_k\subseteq2^{[n]}$ be uncomplemented and pairwise incomparable families and for every pair of families $\mathcal{A}_i$ and $\mathcal{A}_j$ such that $i\not=j$ there exists some $A\in \mathcal{A}_i$ and $\overline{A}\in \mathcal{A}_j$ then it is the case that \begin{equation*}
     \sum^k_{i=1}|\mathcal{A}_i|\leq 2^{n-1}-2
 \end{equation*}
    
\end{theorem}

\section{Preliminaries}

 For a ground set $X$ and a family $\mathcal A\subseteq2^X$, we let $\overline A = X\setminus A$ and $\overline{\mathcal A}=\{\overline A \,:\,A\in\mathcal A\}$. When we will be looking at complements with respect to two distinct sets $X$ and $Y$, we will at times use $\overline{A}$ and $A^c$ respectively.
We also let $\mathcal{U}_{\mathcal{A}}=\{B\subseteq X:A\subseteq B\text{ for some }A\in\mathcal A\}$ be its \emph{upset}, and $\mathcal{D}_{\mathcal{A}}=\{B\subseteq X:B\subseteq A\text{ for some }A\in\mathcal A\}$ be its \emph{downset}. We define the \emph{comparability graph} ${G}_\mathcal{A}$ to have vertex set $Q({G_\mathcal{A}})=\mathcal{A}$ and edge set $E(G_\mathcal{A})=\{AB:A,B\in \mathcal{A}\text{ and } A\subset B \text{ or } B\subset A\}$. Where clear we will refer to the comparability graph $G_\mathcal{A}$ implicitly with its vertex set $\mathcal{A}$.  We let $s(\mathcal{A})$ to be set of complementary pairs $\{A,\overline A \}$ with $A,\overline A\in\mathcal A$ such that $A$ and $ \overline A$ are in the \emph{same} connected component of $G(\mathcal A)$. We let $d({\mathcal{A}})$ to be the set of complementary pairs $\{A,\overline A \}$ with $A,\overline A\in\mathcal A$ such that $A, \overline A$ are contained in \emph{different} connected components of $G(\mathcal A)$. We then therefore have that $|s(\mathcal{A})|+|d({\mathcal{A}})|=\frac{1}{2}|\mathcal A\cap\overline{\mathcal A}|$. For $\mathcal A\subseteq2^X$, we define its \emph{comparability shadow} by $\Gamma(\mathcal A)=\mathcal{U}_\mathcal A\cup \mathcal{D}_\mathcal A$. We then define $b(\mathcal A)=\{\{S,\overline S\}\,:\, S ,\overline S\in\Gamma(\mathcal A)\}$ and note that $ |b(\mathcal A)|$ counts the number of complementary pairs present in the comparability shadow of $\mathcal{A}$. For a family $\mathcal{A}\subseteq 2 ^{[n]}$, and $x\in [n]$ and $Q=[n]\setminus x$ we let $\mathcal{A}_0^{x}=\{A\subseteq Q:A\in\mathcal{A}\}$ and $\mathcal{A}_1^{x}=\{A\subseteq Q:A\cup {x}\in\mathcal{A}\}$. A family $\mathcal A\subseteq2^{[n]}$ is called \emph{convex} whenever for $A,B\in\mathcal A$ and $A\subseteq X\subseteq B$, we have that $X\in\mathcal A$. For an arbitrary family $\mathcal G\subseteq2^{[n]}$, we let its \emph{convex expansion} be $con(\mathcal G)={\mathcal{U}_{\mathcal{G}}}\cap{\mathcal{D}_{\mathcal{G}}}$. So we have that $con(\mathcal G)$ contains all sets $X$ for which there exist $A,B\in\mathcal G$ such that $A\subseteq X\subseteq B$. For $\mathcal{D},\mathcal{I}\subseteq2^{[n]}$  let $\Delta(\mathcal{D},{\mathcal{I}})=|\mathcal{D}\cap\mathcal{I}|-|\mathcal{D}\cap\overline{\mathcal{I}}|-|d(\mathcal{D}\cap\overline{\mathcal{I}})|.$ We let $\gamma_x=|d(C_0^x)|+|d(C_1^x)|+|R_x|-|d(C)|$. We will also at times write  $\Delta(\mathcal{D},\mathcal{I})$ as $\Delta(\mathcal{D}\cap\overline{\mathcal{I}})$ when $\mathcal{D}$ and $\mathcal{I}$ are clear.
In some of the proofs we will use the following theorem of Harris~\cite{Harris1960} and Kleitman~\cite{Kleitman1966}.
\begin{theorem}\label{Harris}
Let $\mathcal{D},\mathcal{I}\subseteq 2^{[n]}$ be two upsets or two downsets. Then it is the case that
\begin{equation*}
    |\mathcal{D}\cap\mathcal{I}|
    \geq
    \frac{|\mathcal{D}|\,|\mathcal{I}|}{2^n},
\end{equation*}
and if $\mathcal{U}\subseteq 2^{[n]}$ is an upset and $\mathcal{D}\subseteq 2^{[n]}$ is a downset then
\begin{equation*}
|\mathcal{U}\cap\mathcal{D}|
    \leq
    \frac{|\mathcal{U}|\,|\mathcal{D}|}{2^n}.
\end{equation*}
\end{theorem}

\section{Proof of Conjecture \ref{conj}}\label{firstsection}
In order to prove Conjecture \ref{conj} we would like to prove that for $\mathcal{G}=\bigcup_{j}\mathcal A_j$ that $|b(\mathcal{G})|\geq |\mathcal{G}|$. Since $2^{n-1}\geq |b(G)|$ we then clearly get that $2^{n-1}\geq |b(G)|\geq |G|$ for which is the conjectured bound. In order to prove that $|b(\mathcal{G})|\geq |G|$, we would like to focus on proving results on the convex expansion $\mathcal{C}=con(\mathcal{G})$, and then show that 
if $|b(\mathcal{C})|\geq |\mathcal{C}|-|s(\mathcal{C})|$ then we also get that $|b(\mathcal{G})|\geq |\mathcal{G}|-|s(\mathcal{G})|$. To do this we define $P=\mathcal{D}_{\mathcal{C}}\cap\overline{\mathcal{U}_{\mathcal{C}}}$. For every $A\in P$ we have $A\in\mathcal{D}_{\mathcal{C}}\subseteq\Gamma(\mathcal{C})$ and $\overline{A}\in\mathcal{U}_{\mathcal{C}}\subseteq\Gamma(\mathcal{C})$, and so $\{A,\overline{A}\}\in b(\mathcal{C})$. This then implies that $|b(\mathcal{C})|\geq |P|-\frac{1}{2}|P\cap \overline{P}|$ and using Lemma \ref{othercomp}, Lemma \ref{down} bounds $|P|$ with respect to $|\mathcal{C}|$ and $|d(\mathcal{{C}})|$.

\begin{lemma}
\label{othercomp}
Let $\mathcal{D},\mathcal{I}\subseteq2^{[n]}$ be downsets. We fix an arbitrary $x\in[n]$, and let $Q=[n]\setminus\{x\}$ and $\mathcal{D}_0^x,\mathcal{D}_1^x,\mathcal{I}_0^x,\mathcal{I}_1^x\subseteq2^Q$. Let $\mathcal C=\mathcal{D}\cap\overline{\mathcal{I}}\subseteq2^{[n]}$ and let $\mathcal C_0^x=\mathcal{D}_0^x\cap\{Q\setminus I:I\in{\mathcal{I}_1^x}\}$, $\mathcal C_1^x=\mathcal{D}_1^x\cap\{Q\setminus I:I\in{\mathcal{I}_0^x}\}$ and $\mathcal R_x=(\mathcal{D}_0^x\setminus\mathcal{D}_1^x)\cap(\mathcal{I}_0^x\setminus\mathcal{I}_1^x)$.
We then get that
\begin{equation*}
    |d(\mathcal{C})|\le |d(\mathcal C_0^x)|+ |d(\mathcal C_1^x)|+|\mathcal R_x|.
\end{equation*}

\end{lemma}

\begin{proof}
Let ${A,\overline A}$ be a complementary pair in $2^Q$. Then it is the case that complementary pairs $\{A, \overline A\cup\{x\}\}$ and $\{ \overline A,A\cup\{x\}\}$ exist in $2^{[n]}$. If such new complementary pairs exist in $\mathcal{C}$ and lie in different connected components of $G(\mathcal{C})$, we will classify how they are counted with respect to $G({\mathcal{C}_0^x})$ and $G({\mathcal{C}_1^x})$. If $n=1$ then $d(\mathcal{C})=\emptyset$, as $\emptyset\subseteq[1]$, so we may assume that $n\geq 2$, and hence $A\neq\overline{A}$ for every $A\subseteq Q$.

For $A,\overline A\cup\{x\}\in \mathcal{C}$ living in different connected components, we then have that  $A\in\mathcal C_0$ and $\overline A\in\mathcal C_1$. If $A\in\mathcal C_1$, then both $A$ and $ \overline A$ lie in $\mathcal C_1$. Assume that $A$ and $\overline A$ lie in the same connected components in $G({\mathcal{C}_1^x})$ , and so there exists some path in $G({\mathcal C_1)}$ from $A$ to $\overline A$. Then adding $x$ to every set of this path gives a path in $G(\mathcal{C})$ from $A\cup\{x\}$ to $\overline{A}\cup\{x\}$, and since $A\subseteq A\cup\{x\}$ with both sets in $\mathcal{C}$, there exists a path from $A$ to $\overline{A}\cup \{x\}$, contradicting the assumption that they lived in different connected components. Hence $\{A,\overline{A}\}\in d(\mathcal{C}_1^x)$.

Similarly, if $ \overline A\in\mathcal C_0$, assume that $A,\overline{A}$ are in the same connected component of $G(\mathcal{C}_0^x)$; then there exists a path from $A$ to $ \overline A$ in $G(\mathcal{C})$, and since $\overline{A}\subseteq\overline{A}\cup\{x\}$ with both sets in $\mathcal{C}$, a path from $A$ to $\overline{A}\cup \{x\}$, contradicting the assumption that they were in different connected components. Hence $\{A,\overline{A}\}\in d(\mathcal{C}_0^x)$.

For the case that $A\not\in\mathcal C_1$ and $ \overline A\not\in\mathcal C_0$, we have that from the fact that $A\in\mathcal C_0$ and $ \overline A\in\mathcal C_1$, we have $A\in\mathcal{D}_0^x$, $ \overline A\in\mathcal{I}_1^x$, $ \overline A\in\mathcal{D}_1^x$, and $A\in\mathcal{I}_0^x$. Since $\mathcal{I}_1^x\subseteq\mathcal{I}_0^x$, the condition $A\not\in\mathcal C_1$ is equivalent here to $A\not\in\mathcal{D}_1^x$. Additionally, since $\mathcal{D}_1^x\subseteq\mathcal{D}_0^x$, the condition $ \overline A\not\in\mathcal C_0$ is equivalent to $A\not\in\mathcal{I}_1^x$.

Therefore we get that, $A\in(\mathcal{D}_0^x\setminus\mathcal{D}_1^x)\cap(\mathcal{I}_0^x\setminus\mathcal{I}_1^x)=\mathcal R_x$. Therefore for any such $A,\overline{A}\cup \{x\}$ that live in different connected components of $G_\mathcal{C}$ we have that $\{A,\overline{A}\}\in d(\mathcal C_0^x)\cup d(\mathcal C_1^x)$ or $A\in\mathcal R_x$.

It remains to check that no term on the right-hand side is used twice. Every complementary pair of $2^{[n]}$ is of the form $\{A,\overline A\cup\{x\}\}$ for exactly one $A\subseteq Q$, so each complementary pair $\{A,\overline A\}$ of $2^Q$ gives rise to exactly two complementary pairs of $2^{[n]}$, namely $\{A,\overline A\cup\{x\}\}$ and $\{\overline A,A\cup\{x\}\}$. By the above, if exactly one of them lies in $d(\mathcal C)$, it is accounted for by $\{A,\overline A\}\in d(\mathcal C_0^x)\cup d(\mathcal C_1^x)$ or by an element of $\mathcal R_x\cap\{A,\overline A\}$. If both lie in $d(\mathcal C)$, then $A,\overline A,A\cup\{x\},\overline A\cup\{x\}\in\mathcal C$, and as $A\cup\{x\}$ lies in the component of $A$ and $\overline A\cup\{x\}$ lies in the component of $\overline A$, the sets $A,\overline A$ lie in different components of $G(\mathcal C)$, and so do $A\cup\{x\},\overline A\cup\{x\}$. As paths in $G(\mathcal C_0^x)$ and in $G(\mathcal C_1^x)$ give paths in $G(\mathcal C)$, we get $\{A,\overline A\}\in d(\mathcal C_0^x)\cap d(\mathcal C_1^x)$, which accounts for both. Summing over all complementary pairs $\{A,\overline A\}$ of $2^Q$ gives us \begin{equation*}
 |d({\mathcal{C}})|\le |d(\mathcal C_0^x)|+ |d(\mathcal C_1^x)|+|\mathcal R_x|.
\end{equation*}
\end{proof}
We note that the injection from Lemma \ref{othercomp} allows us to track exactly where complementary pairs in $\mathcal{C}$ land, and how they are accounted for with respect to $d(\mathcal{C}_0^x),d(\mathcal{C}_1^x)$ and $\mathcal{R}_x$, therefore having $\gamma_x=0$ implies that every complementary pair found in different components of $\mathcal{C}_0^x$ and $\mathcal{C}_1^x$ must be represented as complementary pairs in different components of $\mathcal{C}$, and any set $A\in \mathcal{R}_x$ must also have its complement $\overline{A}\in \mathcal{C}$.
\begin{lemma}
\label{down}
Let $\mathcal{D},\mathcal{I}\subseteq2^{[n]}$ be downsets, and let $\mathcal C=\mathcal{D}\cap\overline{\mathcal{I}}$. Then we get that \begin{equation*}
    |\mathcal{D}\cap\mathcal{I}|\ge|\mathcal C|+|d(\mathcal{C})|.
\end{equation*}
\end{lemma}

\begin{proof}
We do induction on $n$. The case $n=0$ is immediate. Let $n\ge1$, and fix some $x\in[n]$ and let $Q=[n]\setminus x$. We let $\mathcal C_0=\mathcal{D}_0^x\cap\overline{\mathcal{I}_1^x}$ and $\mathcal C_1=\mathcal{D}_1^x\cap\overline{\mathcal{I}_0^x}$, where complements are taken with respect to $Q$. Let $\mathcal R_x=(\mathcal{D}_0^x\setminus\mathcal{D}_1^x)\cap(\mathcal{I}_0^x\setminus\mathcal{I}_1^x)$.
By the induction hypothesis applied in $2^Q$ to the downsets $\mathcal{D}_0^x,\mathcal{I}_1^x$, we have $|\mathcal{D}_0^x\cap\mathcal{I}_1^x|\ge|\mathcal C_0|+ |d(\mathcal C_0)|$and again applying the induction hypothesis to  $\mathcal{D}_1^x,\mathcal{I}_0^x$ gives us $|\mathcal{D}_1^x\cap\mathcal{I}_0^x|\ge|\mathcal C_1|+ |d(\mathcal C_1)|$.
Since $\mathcal{D}_1^x\subseteq\mathcal{D}_0^x$ and $\mathcal{I}_1^x\subseteq\mathcal{I}_0^x$, we then get,
\begin{equation*}
|\mathcal{D}_0^x\cap\mathcal{I}_0^x|+|\mathcal{D}_1^x\cap\mathcal{I}_1^x|
=
|\mathcal{D}_0^x\cap\mathcal{I}_1^x|+|\mathcal{D}_1^x\cap\mathcal{I}_0^x|\
+|(\mathcal{D}_0^x/\mathcal{D}_1^x)\cap(\mathcal{I}_0^x/\mathcal{I}_1^x)|.
\end{equation*}
We have that $|\mathcal{D}_0^x\cap\mathcal{I}_0^x|+|\mathcal{D}_1^x\cap\mathcal{I}_1^x|=$ $|\mathcal{D}\cap\mathcal{I}|$, and $|\mathcal C|=|\mathcal C_0|+|\mathcal C_1|$. We then get that $|D\cap I|\geq |C|+|d(\mathcal{C}_1^x)|+|d(\mathcal{C}_0^x)|+|R_x|$
Then applying Lemma \ref{othercomp} we get that $|\mathcal{D}\cap\mathcal{I}|\ge|\mathcal C|+|d({\mathcal{C}})|$.

\end{proof}
We note that for a convex $\mathcal{C}$ we will write $\Delta(\mathcal{C})=\Delta(D_{\mathcal{C}},\overline{U_{\mathcal{C}}})$ and for a slice we will take $\Delta(\mathcal{C}_0^x)=\Delta(D_{0}^x,I_1^x)$ and $\Delta(\mathcal{C}_1^x)=\Delta(D_{1}^x,I_0^x)$. We also note that for any $D,\overline{I}$ such that $\mathcal{C}=D\cap \overline{I}$ we have that $\Delta(D_{\mathcal{C}},\overline{U_{\mathcal{C}}})\leq \Delta(D,I)$.

For $\mathcal{G}\subseteq 2^{[n]}$ and $\mathcal{C}=con(\mathcal{G})$, \ref{convexblock} gives us a relationship between the sizes of $b(\mathcal{C})$ versus $b(\mathcal{G})$ and $|s(\mathcal{C})|$ versus $s(\mathcal{G})$. We first note that for $\mathcal G\subseteq2^{[n]}$ and $\mathcal{C}=con(\mathcal{G})$ it is the case that $\mathcal{U_C}=\mathcal{U_G}$ and $\mathcal{D_{\mathcal{C}}} =\mathcal {D_G}$. Therefore, $|b(\mathcal C)|=|b(\mathcal G)|$. This is the case as $\mathcal G\subseteq\mathcal C\subseteq\mathcal U_G$, and so we have $\mathcal U_{\mathcal{C}}=\mathcal U_G$. The proof for downward equality is symmetric. Therefore $\Gamma(\mathcal C)=\Gamma(\mathcal G)$, and therefore $|b(\mathcal C)|=|b(\mathcal G)|$.

\begin{lemma}
\label{convexblock}
For $\mathcal G\subseteq2^{[n]}$ and $\mathcal{C}=con(\mathcal{G})$ every connected component of $G(\mathcal{C})$ contains exactly one connected component of $G(\mathcal G)$. Additionally $|s({\mathcal{C}})|-|s({\mathcal{G}})|\le|\mathcal C\setminus\mathcal G|$.
\end{lemma}

\begin{proof}
Let $X\in\mathcal C$. Choose $A,B\in\mathcal G$ with $A\subseteq X\subseteq B$. Then $A$ and $B$ are comparable, so they are in the same connected component of $G(\mathcal G)$. Now suppose $X,Y\in\mathcal C$ are comparable. Without loss of generality, assume $X\subseteq Y$ and let  $A_X,B_X,A_Y,B_Y\in\mathcal G$ be such that $A_X\subseteq X\subseteq B_X$ and $A_Y\subseteq Y\subseteq B_Y$. Then $A_X\subseteq B_Y$, so $X$ and $Y$ are attached to the same component of $G(\mathcal G)$.

Then it is the case that $G(\mathcal C)$ never join sets attached to different components of $\mathcal G$. Additionally any element of $\mathcal C$ is comparable to witnesses in its attached component of $\mathcal G$. Therefore the connected components of $\mathcal C$ are obtained by the convex expansion of the connected components of $\mathcal G$ one by one. 

For a complementary pair in $s({\mathcal{C}})$ but not in $s({\mathcal{G}})$, since convex expansion does not merge connected components, if both sets of this complementary pair belonged to $\mathcal G$, then the pair would already be counted by $s({\mathcal{G}})$. Therefore at least one set in said pair belongs to $\mathcal C\setminus\mathcal G$. Choosing one representative for each such pairs gives an injection into $\mathcal C\setminus\mathcal G$. Therefore we get that 
$|s({\mathcal{C}})|-|s({\mathcal{G}})|\le|\mathcal C\setminus\mathcal G|$.
\end{proof}

\begin{proof}[Proof of Conjecture \ref{conj}]
 Let $\mathcal{G}=\bigcup_j\mathcal{A}_j$ and $\mathcal{C}=con(\mathcal{G})$. Now we know that for $P=\mathcal{D}_{\mathcal{C}}\cap \overline{\mathcal{U}_{\mathcal{C}}}$  that $|b(\mathcal C)|\geq |P|-\frac{1}{2}|P\cap\overline{P}|$, and applying Lemma \ref{down} to the downsets $\mathcal{D}_{\mathcal{C}}$ and $\overline{\mathcal{U}_{\mathcal{C}}}$, for which $\mathcal{D}_{\mathcal{C}}\cap\mathcal{U}_{\mathcal{C}}=\mathcal{D}_{\mathcal{G}}\cap\mathcal{U}_{\mathcal{G}}=\mathcal{C}$, we get that $|P|\geq |\mathcal C|+|d(\mathcal C)|$.
If $S\in\mathcal P\cap\overline{\mathcal P}$, then $S\in\mathcal{D_{\mathcal{C}}}$ and $\overline S\in\mathcal P\subseteq\overline{\mathcal U_{\mathcal{C}}}$ and since $\overline S\in\overline{\mathcal{U_{\mathcal{C}}}}$, we get $S\in\mathcal{U_{\mathcal{C}}}$ and so $S\in\mathcal{D_{\mathcal{C}}}\cap\mathcal{U_{\mathcal{C}}}=\mathcal C$. Similarly, $\overline S\in\mathcal C$. Therefore it is the case that $S\in\mathcal C\cap\overline{\mathcal C}$. The reverse is clear and so we have that $P\cap \overline{P}=\mathcal C\cap\overline{\mathcal C}$, giving $\frac12|P\cap \overline{P}|=\frac12|\mathcal C\cap\overline{\mathcal C}|=|d(\mathcal C)|+|s(\mathcal C)|$. Hence $|b(\mathcal C)|\geq |P|-|d(\mathcal C)|-|s(\mathcal C)|\geq|\mathcal C|-|s(\mathcal C)|$. Therefore, as $|b(\mathcal G)|=|b(\mathcal C)|$, we have by Lemma \ref{convexblock} that $|b(\mathcal G)|\geq |\mathcal G|-|s(\mathcal G)|$. Since the families $\mathcal{A}_i$ are pairwise incomparable, every connected component of $G(\mathcal G)$ lies inside a single $\mathcal A_i$, and since the $\mathcal{A}_i$ are uncomplemented it is the case that $s(\mathcal G)=\emptyset $, and so we get that \begin{equation*}
2^{n-1}\geq |b(\mathcal G)|\geq |\mathcal G|=\sum_{i=1}^k|\mathcal A_i|.
\end{equation*}
\end{proof}
\section{Proof of Theorem \ref{Main}}\label{secondsection}
In order to then show that we get in fact $\sum |\mathcal{A}_i|\leq 2^{n-1}-1$ (proving Theorem \ref{Main})
we first characterise the comparability graph of the intersection of up and down sets that reach equality in the Harris-Kleitman inequality in Theorem \ref{Harris}.
\begin{lemma}
\label{equality Harris}
For $n\in\mathbb{N}$ let $\mathcal{D}\subseteq2^{[n]}$ be a downset and let $\mathcal{U}\subseteq2^{[n]}$ be an upset. If $|\mathcal{D}\cap \mathcal{U}|= \frac{|\mathcal{D}||\mathcal{U}|}{2^{n}}$, then either $G(\mathcal{D}\cap \mathcal{U})$ is connected or $\mathcal{D}\cap \mathcal{U}=\emptyset$.
\end{lemma}

\begin{proof}
We do induction on $n$, and for $n=0$ this is clear. Fix $x\in [n]$, and let $Q=[n]\setminus x$.
As $\mathcal{D}$ is a downset we get that $\mathcal{D}_1^x\subseteq \mathcal{D}_0^x$ and as $\mathcal{U}$ is an upset, $\mathcal{U}_0^x\subseteq \mathcal{U}_1^x$.

By the Harris-Kleitman inequality \ref{Harris} we get that for $i\in \{0,1\}$ that  $|\mathcal{D}_i^x\cap \mathcal{U}_i^x|\le \frac{|\mathcal{D}_i^x||\mathcal{U}_i^x|}{2^{|Q|}}$.

We also have that \begin{equation*}
 |\mathcal{D}\cap \mathcal{U}|=|\mathcal{D}_0^x\cap \mathcal{U}_0^x|+|\mathcal{D}_1^x\cap \mathcal{U}_1^x|\leq\frac{|\mathcal{D}_0^x||\mathcal{U}_0^x|+|\mathcal{D}_1^x||\mathcal{U}_1^x|}{2^{|Q|}}.
\end{equation*}
Since $|\mathcal{U}_1^x|\ge |\mathcal{U}_0^x|$ and $|\mathcal{D}_0^x|\ge |\mathcal{D}_1^x|$, we have $(|\mathcal{D}_0^x|-|\mathcal{D}_1^x|)(|\mathcal{U}_1^x|-|\mathcal{U}_0^x|)\geq0$.
A direct expansion gives \begin{equation*}\frac{|\mathcal{D}_0^x||\mathcal{U}_0^x|+|\mathcal{D}_1^x||\mathcal{U}_1^x|}{2^{|Q|}}=\frac{(|\mathcal{D}_0^x|+|\mathcal{D}_1^x|)(|\mathcal{U}_0^x|+|\mathcal{U}_1^x|)}{2^{|Q|+1}}-\frac{(|\mathcal{D}_0^x|-|\mathcal{D}_1^x|)(|\mathcal{U}_1^x|-|\mathcal{U}_0^x|)}{2^{|Q|+1}},\end{equation*} and the first term on the right is $|\mathcal{D}||\mathcal{U}|/2^{n}=|\mathcal{D}\cap \mathcal{U}|$. Combining this with the previous equation we get that $(|\mathcal{D}_0^x|-|\mathcal{D}_1^x|)(|\mathcal{U}_1^x|-|\mathcal{U}_0^x|)\leq 0$ and so we
have for $i\in \{0,1\}$ that $ |\mathcal{D}_i^x\cap \mathcal{U}_i^x|= \frac{|\mathcal{D}_i^x||\mathcal{U}_i^x|}{2^{|Q|}}$
and $(|\mathcal{D}_0^x|-|\mathcal{D}_1^x|)(|\mathcal{U}_1^x|-|\mathcal{U}_0^x|)=0$. So by induction each of $G(\mathcal{D}_0^x\cap\mathcal{U}_0^x)$ and $G(\mathcal{D}_1^x\cap\mathcal{U}_1^x)$ is connected or empty, and if one of them is empty there is nothing more to prove.
So either $|\mathcal{D}_0^x|=|\mathcal{D}_1^x|$ or $|\mathcal{U}_0^x|=|\mathcal{U}_1^x|$ and so either $\mathcal{D}_0^x=\mathcal{D}_1^x$ or $\mathcal{U}_0^x=\mathcal{U}_1^x$. Assume that $\mathcal{D}_0^x=\mathcal{D}_1^x$, then we have that $\mathcal{D}_0^x\cap \mathcal{U}_0^x=\mathcal{D}_1^x\cap\mathcal{U}_0^x$, and by the fact that $\mathcal{U}$ is an upset, we get that for $A\in \mathcal{D}_1^x\cap\mathcal{U}_0^x$ that $A\in\mathcal{D}_1^x\cap\mathcal{U}_1^x$, so that $A$ and $A\cup\{x\}$ both lie in $\mathcal{D}\cap\mathcal{U}$, connecting the two connected components together. The case where $\mathcal{U}_1^x=\mathcal{U}_0^x$ is similar and so $G(\mathcal{D}\cap\mathcal{U})$ is connected or empty.
\end{proof}
\begin{lemma}
\label{downsetequality}
If for $D,I\subseteq 2^{[n]}$ being downsets it is the case that $|D\cap I |=|D||I|/2^{n}$, then $\{x:|D_0^x|>|D_1^x|\}\cap\{x:|I_0^x|>|I_1^x|\}=\emptyset$
\end{lemma}
\begin{proof}
Since both $D$ and $I$ are downsets we have that $|D_0^x|\geq|D_1^x|$ and $|I_0^x|\geq |I_1^x|$ and so, applying Theorem \ref{Harris} in each of the two slices,\begin{equation*}
|D\cap I|\geq \frac{|D_0^x||I_0^x|+|D_1^x||I_1^x|}{2^{n-1}}=\frac{|D||I|}{2^n}+\frac{(|D_0^x|-|D_1^x|)(|I_0^x|-|I_1^x|)}{2^n}
\end{equation*}
Now since we have that $|D\cap I |=|D||I|/2^{n}$ we have that $(|D_0^x|-|D_1^x|)(|I_0^x|-|I_1^x|)=0$, and so we have that $\{x:|D_0^x|>|D_1^x|\}\cap\{x:|I_0^x|>|I_1^x|\}=\emptyset$
\end{proof}
As mentioned earlier, we aim to show that by fixing $\mathcal{A}_1$ and letting $\mathcal{G}=\bigcup_{i\not=1}\mathcal{A}_i$ and $\mathcal{C}=con({\mathcal{G}})=D_{\mathcal{G}}\cap U_{\mathcal{G}}$ that we get $|D_{\mathcal{G}}\cap \overline{U_{\mathcal{G}}}|\geq |\mathcal{C}|+|d(\mathcal{C})|+1$. To do so we aim at showing in Lemma \ref{downsetint}, by induction, that for general $\mathcal{D},\mathcal{I}\subseteq2^{[n]}$ downsets, and $\mathcal C=\mathcal{D}\cap\overline{ \mathcal{I}}$ that if $G(\mathcal C)$ is disconnected and $\mathcal{D}\cup\overline{ \mathcal{I}}\not=2^{[n]}$, then
$|\mathcal{D}\cap \mathcal{I}|\ge|\mathcal C|+|d(\mathcal C)|+1$. In preparation for Lemma \ref{downsetint}, we prove Lemmas \ref{Induction step} and \ref{bound} that will be used in the proof.
\begin{lemma}
\label{Induction step}
 For $n\in\mathbb{N}$, $Q\subset[n]$, let $\mathcal L\subseteq2^Q$ be non-empty and convex with $G(\mathcal{L})$ connected, and let $X\subseteq Q$ be incomparable with every element of $\mathcal L$. Let
$\mathcal S=\mathcal{D}_{\mathcal L\cup\{X\}}\cap\overline{\mathcal{U}_{\mathcal L\cup\{X\}}}$, where complements are taken with respect to $Q$, then it is the case that $|\mathcal S|\ge|\mathcal L|+2$.
\end{lemma}

\begin{proof}
Since $X$ is incomparable with every element of $\mathcal L$, the graph $G(\mathcal L\cup\{X\})$ is disconnected and $\mathcal{L}\cup\{X\}$ is convex. Let $D=\mathcal{D}_{\mathcal{L}\cup \{X\}}$ and $U=\mathcal{U}_{\mathcal L\cup \{X\}}$, and so $\mathcal S=D\cap\overline U$.
Since $D$ and $\overline U$ are both downsets, the Harris-Kleitman inequality \ref{Harris} gives \begin{equation*}
     |D\cap\overline U|\ge\frac{|D||\overline U|}{2^{|Q|}}=\frac{|D||U|}{2^{|Q|}}.
\end{equation*} As $\mathcal L\cup\{X\}$ is convex we have $D\cap U=\mathcal L\cup\{X\}$, which is non-empty and disconnected, so by Theorem \ref{Harris} and Lemma \ref{equality Harris} (applied in $2^Q$) we get that \begin{equation*}|D\cap U|<\frac{|D||U|}{2^{|Q|}}\quad\end{equation*}
and so $|D||U|>2^{|Q|}|\mathcal L\cup \{X\}|$.
Therefore $|\mathcal S|=|D\cap\overline U|> |\mathcal L|+1$ and so $|\mathcal S|\ge|\mathcal L|+2$.
\end{proof}

\begin{lemma}
\label{bound}
Let $\mathcal C\subseteq2^{[n]}$ be convex and $P=\mathcal{D}_{\mathcal C}\cap\overline{\mathcal{U}_{\mathcal C}}$.
Let $G(\mathcal C)$ have $r\geq 2$ connected components $T_1,..., T_r$, and $W\not\in\Gamma(\mathcal C)$. Let there be disjoint sets \begin{equation*}
M_i=\{x\not\in W:\exists A\in T_i,\ x\not\in A\}\cup\{x\in W:\exists A\in T_i,\ x\in A\}.
\end{equation*} and $\overline{W}_i=\overline{W}\cap M_i$ and $W_i=M_i\setminus \overline{W}_i$.
Let $\mathcal L_i=\{A\cap M_i:A\in T_i\}$ and $\mathcal S_i=\mathcal{D}_{\mathcal L_i\cup\{\overline{W}_i\}}\cap\{M_i\setminus Y:Y\in \mathcal{U}_{\mathcal L_i\cup\{\overline{W}_i\}}\}$, with the downset and upset taken in $2^{M_i}$. Then $T_i=\{L\cup (\overline{W}\setminus M_i):L\in\mathcal L_i\}$ and $\mathcal{S}_i\subseteq P$.
\end{lemma}

\begin{proof}
 By definition of $M_i$, every $A\in T_i$ satisfies $A\setminus M_i=\overline{W}\setminus M_i$, so $T_i=\{(\overline{W}\setminus M_i)\cup L:L\in\mathcal L_i\}$,
where $\mathcal L_i\subseteq2^{M_i}$ is convex and connected. Recall that
$$\mathcal S_i=
\mathcal{D}_{\mathcal L_i\cup\{\overline{W}_i\}}
\cap
\left\{M_i\setminus Y:Y\in \mathcal{U}_{\mathcal L_i\cup\{\overline{W}_i\}}\right\}.
$$
For $S\in \mathcal S_i$, there are some $X,Y\in\mathcal L_i\cup\{\overline{W}_i\}$ such that $S\subseteq X
$ and $
S\cap Y=\emptyset$. If $X,Y\in\mathcal L_i$, then there are 
$
A=(\overline{W}\setminus M_i)\cup X\in T_i
$ and $
B=(\overline{W}\setminus M_i)\cup Y\in T_i$
. As we have that $S\subseteq M_i$, the fact that $S\subseteq X$ and $S\cap Y=\emptyset$ implies that
$S\subseteq A$ and $S\cap B=\emptyset$ and so $S\in P$. If it is the case that $X\in\mathcal L_i$ and $Y=\overline{W}_i$. Then there exists some $A=(\overline{W}\setminus M_i)\cup X\in T_i$ and so $S\subseteq A$. Since $S\cap \overline{W}_i=\emptyset$, we can take any component $T_j$ with $j\not= i$ and any $B=\overline{W}\setminus M_j\cup L\in T_j$ and $L\in \mathcal{L}_j$. As $M_i\cap M_j=\emptyset$, every element $K\in T_j$ has the property that for all $x\not\in M_j$ then $x\in K$ if and only if $x\in\overline{W}$. Therefore we have that $B\cap M_i=\overline{W}_i$. As $S\subseteq M_i$, $L\subseteq M_j$, $M_i\cap M_j=\emptyset$ and $S\cap \overline{W}_i=\emptyset$, we get that $
S\cap B=\emptyset$ and so $S\in P$. For $X=\overline{W}_i$ and $Y\in\mathcal L_i$, we get that $S\subseteq \overline{W}_i$ and so for $j\not=i$ and $B=\overline{W}\setminus M_j\cup L\in T_j$ we get that $S\subseteq B$. Since $S\cap Y=\emptyset$ and $S\subseteq M_i$ we have for $A=\overline{W}\setminus M_i\cup Y\in T_i$ that $S\cap A=\emptyset$ and so $S\in P$.
If $X=Y=\overline{W}_i$, then $S\subseteq \overline{W}_i$ and $S\cap \overline{W}_i=\emptyset$, so $S=\emptyset$. As $\mathcal C$ is non-empty we get that $\emptyset\in \mathcal{D}_{\mathcal C}$ and $\emptyset\in\overline{\mathcal{U}_{\mathcal C}}$, so $\emptyset\in P$. Therefore we get that for every $S\in\mathcal S_i$ that $S\in P$ and so $\mathcal S_i\subseteq P$.
\end{proof} 

We are now in a position to prove Lemma \ref{downsetint}.
\begin{lemma}
\label{downsetint}
Let $\mathcal{D},\mathcal{I}\subseteq2^{[n]}$ be downsets, and let $\mathcal C=\mathcal{D}\cap\overline{ \mathcal{I}}$.
If $G(\mathcal C)$ is disconnected and $\mathcal{D}\cup\overline{ \mathcal{I}}\not=2^{[n]}$, then
$|\mathcal{D}\cap \mathcal{I}|\ge|\mathcal C|+|d(\mathcal C)|+1$.
\end{lemma}

\begin{proof}
We prove the result by induction on $n$. The case $n=0$ is clear, as then $|\mathcal C|\leq 1$. Assume that for minimal $1\leq n\in \mathbb{N}$ there exists $\mathcal{D},\mathcal{I}\subseteq 2^{[n]}$ such that $\mathcal{D}\cap \mathcal{I}$ is a counterexample. 
By Lemma \ref{down} we get that $\Delta(\mathcal{D},\mathcal{I})\ge0$ and so by $\mathcal{D}\cap\mathcal{I}$ being a counterexample it is the case that $\Delta(\mathcal{D},\mathcal{I})=0$.
Let $W\not\in \mathcal{D}\cup\overline{ \mathcal{I}}$, and take $x\in[n]$, and let $Q=[n]\setminus\{x\}$. As in Lemma \ref{othercomp} we have that
$\mathcal C_0^x=\mathcal{D}_0^x\cap\{Q\setminus I:I\in {\mathcal{I}_1^x}\}$ and
$\mathcal C_1^x=\mathcal{D}_1^x\cap\{Q\setminus I:I\in {\mathcal{I}_0^x}\}$,
and let $\mathcal R_x=(\mathcal{D}_0^x\setminus \mathcal{D}_1^x)\cap(\mathcal{I}_0^x\setminus \mathcal{I}_1^x)$.
By the identity in the proof of Lemma \ref{down}, we then get that
\begin{equation*}
\Delta(\mathcal{D},\mathcal{I})
=
\Delta(\mathcal{D}_0^x,\mathcal{I}_1^x)+\Delta(\mathcal{D}_1^x, \mathcal{I}_0^x)
+
\gamma_x.
\end{equation*}
By Lemma \ref{othercomp}, we get that $
\gamma_x\geq 0
$, and by Lemma \ref{down} the other two terms are non-negative, and so it must be the case that $\Delta(\mathcal{D}_0^x,\mathcal{I}_1^x)=0$ and $\Delta(\mathcal{D}_1^x,\mathcal{I}_0^x)=0$. If $x\not\in W$, then $W\subseteq Q$, $W\not\in \mathcal{D}_0^x$, and $Q\setminus{W}\not\in \mathcal{I}_1^x$. Therefore $W\not\in \mathcal{D}_0^x\cup\{Q\setminus I:I\in\mathcal{I}_1^x\}$.
If $G(\mathcal C_0^x)$ is disconnected, then by induction we would get that $\Delta(\mathcal{D}_0^x,\mathcal{I}_1^x)\ge1$ which is clearly not possible and so $G(\mathcal C_0^x)$ is connected or empty. Similarly, if $x\in W$, then $W\setminus \{x\}\not\in \mathcal{D}_1^x\cup \{Q\setminus M:M\in\mathcal{I}_0^x\}$, and so $G(\mathcal C_1^x)$ is connected or empty. In other words, for every $x\in[n]$ the sets of $\mathcal C$ that agree with $W$ on $x$ induce a connected or empty subgraph of $G(\mathcal C)$.

Let $T_1,\ldots,T_r$ be the connected components of $G(\mathcal C)$, with $r\ge2$. For each component $T_i$ let
\begin{equation*}
M_i=\{x\not\in W:\exists A\in T_i,\ x\not\in A\}\cup\{x\in W:\exists A\in T_i,\ x\in A\}.
\end{equation*}
And so $\overline{M_i}$ will be all the elements of $x\in [n]$ for which for any $A\in T_i$ it is the case that $x\in A$ if and only if $x\in \overline{W}$.

Assume that $M_i\cap M_j\not=\emptyset$ for some $i\neq j$, and let $y\in M_i\cap M_j$. If $y\not\in W$, there are $K_1\in T_i$ and $K_2\in T_j$ with $y\not\in K_1,K_2$, so that $K_1,K_2\in\mathcal C_0^y$. By the above $G(\mathcal C_0^y)$ is connected, and a path in $G(\mathcal C_0^y)$ is a path in $G(\mathcal C)$, contradicting $T_i\neq T_j$. The case when $y\in W$ is similar, using $\mathcal C_1^y$, and so $M_i\cap M_j=\emptyset$.

Let $\overline{W}_i=\overline{W}\cap M_i$ and $\mathcal L_i=\{A\cap M_i:A\in T_i\}$. Then by Lemma \ref{bound} we get that $T_i=\{(\overline{W}\setminus M_i)\cup L:L\in\mathcal L_i\}$,
where $\mathcal L_i\subseteq2^{M_i}$ is convex and connected (a connected component of a convex family is convex) and $|\mathcal L_i|=|T_i|$. Let $P=\mathcal{D}_{\mathcal C}\cap\overline{\mathcal{U}_{\mathcal C}}$ and
since $\mathcal{D}_{\mathcal C}\subseteq \mathcal{D}$ and $\overline{\mathcal{U}_{\mathcal C}}\subseteq \mathcal{I}$, we have $P\subseteq \mathcal{D}\cap \mathcal{I}$.
We then define $$\mathcal S_i=
\mathcal{D}_{\mathcal L_i\cup\{\overline{W}_i\}}
\cap
\left\{M_i\setminus Y:Y\in \mathcal{U}_{\mathcal L_i\cup\{\overline{W}_i\}}\right\}$$ and from Lemma \ref{bound} (note that $\Gamma(\mathcal C)\subseteq\mathcal D\cup\overline{\mathcal I}$, so $W\not\in\Gamma(\mathcal C)$) we get that $\mathcal{S}_i\subseteq P$.

By Lemma \ref{down}, applied in $2^{M_i}$ to the downsets $\mathcal{D}_{\mathcal F_i}$ and $\{M_i\setminus Y:Y\in\mathcal U_{\mathcal F_i}\}$ where $\mathcal F_i=\mathcal L_i\cup\{\overline{W}_i\}$, we get that
$|\mathcal S_i|\ge|\mathcal{D}_{\mathcal F_i}\cap\mathcal U_{\mathcal F_i}|\geq|\mathcal F_i|=|\mathcal L_i|+ \mathbbm{1}[\overline{W}_i\not\in\mathcal L_i]$.
Since the $M_i$ are pairwise disjoint and $\emptyset\in\mathcal S_i\subseteq 2^{M_i}$ for each $i$, the families $\mathcal S_i$ pairwise intersect exactly in $\{\emptyset\}$.

If $|d(\mathcal C)|=0$ and $\overline{W}\not\in\mathcal C$, then it is the case that $\overline{W}_i\not\in\mathcal L_i$ (as $(\overline{W}\setminus M_i)\cup\overline{W}_i=\overline{W}$) and so
\begin{equation*}
|P|\geq |\bigcup_{i=1}^r\mathcal S_i|=1+\sum^r_{i=1}(|\mathcal S_i|-1)\ge1+\sum^r_{i=1}|\mathcal L_i|=|\mathcal C|+1.
\end{equation*}
If $|d(\mathcal C)|=0$ and $\overline{W}\in\mathcal C$, then for some component $T_h$ it is the case that $\overline{W}\in T_h$. For every $i\not= h$ it is the case that $\overline{W}_i$ is incomparable with every element of $\mathcal L_i$, since if $\overline{W}_i$ were comparable with some $L\in\mathcal L_i$ then $\overline{W}$ would be comparable with $(\overline{W}\setminus M_i)\cup L\in T_i$, giving $\overline{W}\in T_i$. Therefore Lemma \ref{Induction step}, applied in $2^{M_i}$ with $\mathcal L=\mathcal L_i$ and $X=\overline{W}_i$, gives us that
$|\mathcal S_i|\ge|\mathcal L_i|+2$ for $i\not= h$.
Since $\overline{W}\in T_h$ we have that $\overline{W}_h\in\mathcal{L}_h$ and so by the bound above $|\mathcal S_h|\geq |\mathcal{L}_h|$. As $r\ge2$ it is the case that
\begin{equation*}
|P|\geq |\bigcup^r_{i=1}\mathcal S_i|
=
1+\sum^r_{i=1}(|\mathcal S_i|-1)
\ge
|\mathcal C|+(r-1)
\ge
|\mathcal C|+1.
\end{equation*}
Finally suppose that $|d(\mathcal C)|\not=0$. Choose a different-component complementary pair $A\in T_i$ and $\overline A\in T_j$. If there is some $x\in[n]$ with $x\not\in M_i\cup M_j$, then $x\in A$ if and only if $x\in \overline{W}$ if and only if $x\in\overline A$, which is not possible, and so $M_i\cup M_j=[n]$.
Since the $M_\ell$ are pairwise disjoint, $M_\ell=\emptyset$ and so $T_\ell=\{\overline{W}\}$ for every $\ell\not\in\{i,j\}$; in particular $r\leq 3$. A component $T_\ell=\{\overline{W}\}$ contains no set whose complement lies in $\mathcal C$, as $W\not\in\mathcal C$. Hence $M_i,M_j\neq\emptyset$, and every different-component complementary pair of $\mathcal C$ lies between $T_i$ and $T_j$. Writing $A=(\overline{W}\setminus M_i)\cup L_1$ and $\overline{A}=(\overline{W}\setminus M_j)\cup L_2$ with $L_1\in\mathcal L_i$ and $L_2\in\mathcal L_j$, complementarity forces $L_1=W\cap M_i$ and $L_2=W\cap M_j$, so the pair is unique and $|d(\mathcal C)|=1$.
Let $W_\ell=W\cap M_\ell$, so that $A=W_i\cup\overline{W}_j$ and $\overline A=\overline{W}_i\cup W_j$. Each of $W_i,\overline{W}_i,W_j,\overline{W}_j$ is non-empty: for instance, if $W_i=\emptyset$ then $W=W_j\subseteq\overline A\in\mathcal D$, and if $\overline{W}_j=\emptyset$ then $A=W_i\subseteq W$, so $W\in\overline{\mathcal I}$, both contradicting the choice of $W$; the other two cases are symmetric. Hence $A$ and $\overline A$ both intersect $M_i$ and $M_j$.

If $r=2$, relabel so that $\{i,j\}=\{1,2\}$. We have seen that
$|\mathcal S_i|\ge|\mathcal L_i|+ \mathbbm{1}[{\overline{W}_i\not\in\mathcal L_i}]$. It is not the case that both $\overline{W}_i\in\mathcal{L}_i$, as then $\overline{W}=\overline{W}_1\cup \overline{W}_2$ would be in both components $T_1$ and $T_2$, contradicting the incomparability condition. We then get that, $
 \mathbbm{1}[\overline{W}_1\not\in\mathcal L_1]+
 \mathbbm{1}[\overline{W}_2\not\in\mathcal L_2]\ge1.
$
As $M_1\cap M_2=\emptyset$, the families $\mathcal S_1$ and $\mathcal S_2$ intersect only in $\emptyset$ and so 
$
|\mathcal S_1\cup\mathcal S_2|=|\mathcal S_1|+|\mathcal S_2|-1\ge
|\mathcal L_1|+|\mathcal L_2|=|\mathcal C|.
$As every element of $\mathcal S_1\cup\mathcal S_2$ is a subset of  $M_1$ or of $M_2$, while $A=W_1\cup \overline{W}_2$ and $\overline A=\overline{W}_1\cup W_2$ both intersect $M_1$ and $M_2$ and so we get that $A,\overline{A}\not\in\mathcal S_1\cup\mathcal S_2$. As $A,\overline{A}\in\mathcal C$ we have $A,\overline{A}\in\mathcal{D}_{\mathcal C}\cap\overline{\mathcal{U}_{\mathcal C}}=P$, and then
$|P|\ge|\mathcal S_1\cup\mathcal S_2|+2\ge|\mathcal C|+2$.

If $r\ge3$, then by the above $r=3$ and the third component is $T_\ell=\{\overline{W}\}$, so that $\overline{W}\in\mathcal C$ and $\mathcal S_\ell=\{\emptyset\}$. Exactly as in the case $\overline{W}\in\mathcal C$ above, Lemma \ref{Induction step} gives $|\mathcal S_i|\geq|\mathcal L_i|+2$ and $|\mathcal S_j|\geq|\mathcal L_j|+2$, and so
$|P|\ge|\bigcup_l\mathcal S_l|=1+(|\mathcal S_i|-1)+(|\mathcal S_j|-1)\geq|\mathcal L_i|+|\mathcal L_j|+3=|\mathcal C|+2=|\mathcal C|+|d(\mathcal C)|+1$.

Therefore in all cases we get that,
$|P|\ge|\mathcal C|+|d(\mathcal C)|+1$ and as 
 $P\subseteq \mathcal{D}\cap \mathcal{I}$, we get that 
$|\mathcal{D}\cap \mathcal{I}|\ge|\mathcal C|+|d(\mathcal C)|+1$,
contradicting the assumption that $\Delta(\mathcal{D},\mathcal{I})=0$.
\end{proof}

\begin{lemma}
\label{finallemma}
Let $\mathcal C\subseteq2^{[n]}$ be convex. If $G(\mathcal C)$ is disconnected and $\Gamma(\mathcal C)\not=2^{[n]}$, then
$|b(\mathcal C)|\ge|\mathcal C|-|s(\mathcal C)|+1$.
\end{lemma}

\begin{proof}
We apply Lemma \ref{downsetint} with
$D=D_{\mathcal C}$ and $\mathcal{I}=\overline{\mathcal{U}_{\mathcal C}}$.
As $\mathcal C$ is convex,
$D\cap\overline{ \mathcal{I}}=D_{\mathcal C}\cap \mathcal{U}_{\mathcal C}=\mathcal C$,
and 
$D\cup\overline{ \mathcal{I}}=D_{\mathcal C}\cup \mathcal{U}_{\mathcal C}=\Gamma(\mathcal C)$.
Letting $P=D_{\mathcal C}\cap\overline{\mathcal{U}_{\mathcal C}}$, we then get that,
$|P|\ge|\mathcal C|+|d(\mathcal C)|+1$.
For $S\in P$, as $S\in\mathcal{D_{\mathcal{C}}}$, it is the case that $S$ is comparable with some element of $\mathcal C$. Since $S\in\overline{\mathcal{U_{\mathcal{C}}}}$, so $\overline S\in\mathcal{U_{\mathcal{C}}}$, and therefore $\overline S$ is also comparable with some element of $\mathcal C$. This then means that every complementary pair represented by a element of $\mathcal P$ is  in $b(\mathcal C)$. The number of distinct complementary pairs indexed by $\mathcal P$ is $|\mathcal P|-\frac12|\mathcal P\cap\overline{\mathcal P}|$.
If $S\in\mathcal P\cap\overline{\mathcal P}$, then $S\in\mathcal{D_{\mathcal{C}}}$ and $\overline S\in\mathcal P\subseteq\overline{\mathcal U_{\mathcal{C}}}$. Since $\overline S\in\overline{\mathcal{U_{\mathcal{C}}}}$, we get $S\in\mathcal{U_{\mathcal{C}}}$. Therefore $S\in\mathcal{D_{\mathcal{C}}}\cap\mathcal{U_{\mathcal{C}}}=\mathcal C$. Similarly, $\overline S\in\mathcal C$. Therefore it is the case that $S\in\mathcal C\cap\overline{\mathcal C}$. The reverse is clear. Therefore we get that $\mathcal P\cap\overline{\mathcal P}=\mathcal C\cap\overline{\mathcal C}.$
For every $S\in P$ we have that $S,\overline S\in \Gamma(\mathcal C)$ and the number of distinct complementary pairs in $P$ is at least
$|P|-|s(\mathcal C)|-|d(\mathcal C)|$ and so
\begin{equation*}
|b(\mathcal C)|
\ge
|P|-|s(\mathcal C)|-|d(\mathcal C)|
\ge
|\mathcal C|-|s(\mathcal C)|+1.
\end{equation*}
\end{proof}

\begin{proof}[Proof of Theorem \ref{Main}]
Relabelling, we may assume that $\mathcal{A}_1$ is non-empty and that at least two of $\mathcal A_2,\dots,\mathcal A_k$ are non-empty. Let $\mathcal G=\bigcup_{j\not=1}\mathcal A_j$ and $ 
\mathcal C=con(\mathcal G)
=
D_{\mathcal G}\cap \mathcal{U}_{\mathcal G}.$
We have that $\Gamma(\mathcal C)=\Gamma(\mathcal G)$ and as every $A\in\mathcal A_1$ is incomparable with $B\in \mathcal G$ it is the case that for $W\in\mathcal A_1$ we have that
$W\not\in\Gamma(\mathcal G)=\Gamma(\mathcal C)$ 
and so $\Gamma(\mathcal C)\not=2^{[n]}$.

As at least two of the pairwise incomparable families $\mathcal A_2,\dots,\mathcal A_k$ are non-empty, the graph $G(\mathcal G)$ has at least 2 connected components and so $G(\mathcal G)$ is disconnected. If $X\in\mathcal{C}$, then $X$ lies between two comparable elements of $\mathcal G$, which must belong to the same component of $G(\mathcal G)$ and two convexified elements assigned to different components cannot be comparable, since that would create a path between the original components. Therefore $G(\mathcal C)$ is also disconnected.

By Lemma \ref{finallemma} it is the case that
$|b(\mathcal C)|\ge|\mathcal C|-|s(\mathcal C)|+1$.
We then have that by Lemma \ref{convexblock} we get that $|s(\mathcal C)|\le|\mathcal C|-|\mathcal G|+|s(\mathcal{G})|$ and so $|b(\mathcal G)|=|b(\mathcal C)|\ge|\mathcal G|-s(\mathcal{G})+1$. Now since $\mathcal{A}_1,...,\mathcal{A}_k$ are all uncomplemented, so are any of the connected components of $G(\mathcal{G})$, and so $|s(\mathcal{G})|=0$. Every complementary pair in $b(\mathcal G)$ has both elements in $\Gamma(\mathcal G)$, and since $\mathcal A_1\cap\Gamma(\mathcal G)=\emptyset$ and $\mathcal A  _1$ is uncomplemented, it can contain at most one set from each of the remaining complementary pairs. Therefore $|\mathcal A _1|\le2^{n-1}-|b(\mathcal G)|$.
And so we get that
\begin{equation*}
\sum_{i=1}^k|\mathcal{A}_i|
=
|\mathcal A_1|+|\mathcal G|
\le
2^{n-1}-|b(\mathcal G)|+|\mathcal G|
\le
2^{n-1}-1.
\end{equation*}
\end{proof}
To conclude this section we provide an example of the tightness of this bound.
\begin{proof}[Proof of Proposition \ref{tight}]
Let $k\ge3$ and take distinct elements $x,x_1,\ldots,x_{k-1}\in[n]$, and let $Q=[n]\setminus\{x\}$. Then for $1\leq i\leq k-1$ we let $\mathcal A_i=\big\{\{x,x_i\}\big\}$ and
$\mathcal A_k=
2^Q\setminus(\{\emptyset\}\cup\bigcup_{1}^{k-1}\{\{{x_i}\}\}).
$
Then all $\mathcal{A}_i$ are non-empty and it is clear that each family is uncomplemented, and pairwise incomparable and so we get that
\begin{equation*}
|\mathcal A_1|+...+|\mathcal A_k|
= (k-1)+(2^{n-1}-1-(k-1))=
2^{n-1}-1.
\end{equation*}
\end{proof}
\section{Proof of Theorem \ref{3rd}}\label{lastsection}
We now require that our families $\mathcal{A}_1,...,\mathcal{A}_k$ must have the property that between any two pairs $\mathcal{A}_i,\mathcal{A}_j$ such that $i\not=j$ it is the case that there exists some $A\in \mathcal{A}_i$ such that $\overline{A}\in\mathcal{A}_j$. Therefore this section is devoted to proving Theorem \ref{3rd}. We now define $\mathcal{G}=\bigcup_{i}\mathcal{A}_i$ and $\mathcal{C}=con(G)$.
In the aim to arrive at a contradiction we assume in this section that $|G|=2^{n-1}-1$. 
\begin{lemma}
Assume that $|\mathcal{G}|=2^{n-1}-1$, then $1\geq\Delta(\mathcal{C})\geq0$.
\end{lemma}
\begin{proof}
For $P=D_{\mathcal{C}}\cap \overline{U_{\mathcal{C}}}$ we have that as before $|b(\mathcal{G})|=|b(\mathcal{C})|\geq |P|-|s(\mathcal{C})|-|d(\mathcal{C})|$ and so $|b(\mathcal{C})|\geq |\mathcal{C}|+\Delta(\mathcal{C})-|s(\mathcal{C})|$. Since $|s(\mathcal{C})|\leq |\mathcal{C}|-|\mathcal{G}|+|s(\mathcal{G})|$ and $s(\mathcal{G})=\emptyset$ we have that $|b(\mathcal{C})|\geq |\mathcal{G}|+\Delta(\mathcal{C})$. Therefore we have that $2^{n-1}\geq |b(\mathcal{C})|\geq 2^{n-1}-1+\Delta(\mathcal{C})$ and so $1\geq \Delta(\mathcal{C})\geq 0$
\end{proof}
\begin{lemma}
\label{fullshadow}
If $G(\mathcal{C})$ is disconnected, and $\Delta(\mathcal{C})=0$ then it is the case that $\Gamma(\mathcal{C})=2^{[n]}$ and $|d(\mathcal{C})|>0$. 
\end{lemma}
\begin{proof}
Since we have by Lemma \ref{downsetint} that if $G(\mathcal{C})$ is disconnected and $\Gamma(\mathcal{C})\not=2^{[n]}$ then we would have that $\Delta(\mathcal{C})\geq 1$ which is a contradiction.

Since $\Delta(\mathcal{C})=0$, we have that $|D_{\mathcal{C}}\cap \overline{U_{\mathcal{C}}}|=|\mathcal{C}|+|d(\mathcal{C})|$, and if $|d(\mathcal{C})|=0$ by the Harris-Kleitman Lemma \ref{Harris} we would get that $|\mathcal{C}|=|D_{\mathcal{C}}\cap \overline{U_{\mathcal{C}}}|\geq |D_{\mathcal{C}}||U_{\mathcal{C}}|/2^{n}\geq |D_{\mathcal{C}}\cap U_{\mathcal{C}}|=|\mathcal{C}|$, and so by Lemma \ref{equality Harris} we would have that $G(\mathcal{C})$ is connected which contradicts the assumption that it was not.
\end{proof}
\begin{lemma}
\label{Variationproperties}
    Given $\mathcal{A}_1,\mathcal{A}_2,\mathcal{A}_3$ that are pairwise incomparable and $G=\bigcup \mathcal{A}_i$, $C=con(G)$ and for every $x\in [n]$ it is the case that $\mathcal{C}_0^x$ and $\mathcal{C}_1^x$ meet $2$ of the families $\mathcal{A}_i$. Then we let $M_i=\{x:\exists A_1,A_2\in \mathcal{A}_i\text{ such that } x\in A_1\text{ and } x\not\in A_2\}$ and when $M_i\not=\emptyset$ we have that for all $j\not=i$ there is some $A\in \mathcal{A}_i$ where $\overline{A}\in\mathcal{A}_j$. Then it is the case that $M_i$ partition $[n]$, are pairwise disjoint and non-empty. Additionally for $A\in A_j$ and $i\not=j$ we let $u_{i,j}=A\cap M_i$, then it is the case that $\emptyset\subsetneq u_{i,j}\subsetneq M_i$.
\end{lemma}
\begin{proof}
We have that for any original family $[n]=\bigcup_i M_i$ since if not then there is some $x\in [n]$ for which  the number of families of $\mathcal{C}_0^x$ and $\mathcal{C}_1^x$ is less than $4$. Now all $M_i$ are pairwise disjoint since then we would have that number of families met by $\mathcal{C}_0^x$ and $\mathcal{C}_1^x$ is $\geq 5$. And so $M_i$ partitions $[n]$ and can be seen as the set of elements of $[n]$ for which the $i$th family \emph{varies} in. 

Additionally all $M_i$ are non-empty, this can be seen as follows; Clearly if $M_i=\emptyset$ then we have that $\mathcal{A}_i$ is a singleton, and so if say $M_1$ was the only set that was empty we would need $2$ distinct sets (coming from the complementary pairs coming from $\mathcal{A}_2$ and $\mathcal{A}_3$) in $\mathcal{A}_1$ which is not possible. If additionally we have that $M_2=\emptyset$ then we would get $\mathcal{A}_1=\{A\}$, $\mathcal{A}_2=\{B\}$ and $M_3=[n]$. As each slice meets $2$ original families, it must be that $A=\overline{B}$. But as $M_3\not=\emptyset$ it is the case that there exists a complementary pair between $\mathcal{A}_3$ and $\mathcal{A}_j$ for $j\in\{1,2\}$ which is not possible.

 Since some family $\mathcal{A}_i$ meets both $\mathcal{C}_0^x$ and $\mathcal{C}_1^x$ for every $x$ we have that it can't be that all $M_i$ are empty. Now since all $M_i$ are non-empty, this then shows that between each pair of families we have atleast some complementary pair in $C$.

We first note that $u_{i,j}=M_i\setminus u_{i,k}$. Assume that without loss of generality $u_{1,2}=\emptyset$ and so $u_{1,3}=M_1$, we then take complementary pairs $A,B\in \mathcal{A}_1,{A^c}\in \mathcal{A}_2$ and ${B^c}\in \mathcal{A}_3$, and so $M_1\subseteq A$ and $B\cap M_1=\emptyset$. We have that for every $A_i,A_j\in \mathcal{A}_1$ that $Z=A_i\cap [n]\setminus M_1=A_j\cap [n]\setminus M_1$ and so $A=Z\cup M_1$ and $B=Z$, and so $B\subseteq A$, and so ${A^c}\subseteq {B^c}$, and since ${A^c}$ and ${B^c}$ lived in different connected components this is a contradiction.
\end{proof}
\begin{lemma}
\label{2comp}
If $\Delta(\mathcal{C})=0$ then there are at most $2$ connected components.
\end{lemma}
\begin{proof}
We will do induction on $n$ where for $n=1$ this is clear. ASsume that $\mathcal{C}$ has at least $3$ components.

Since we have that $\Delta(\mathcal{C})=0$ we have that by Lemma \ref{fullshadow} it is the case that $\Gamma(\mathcal{C})=2^{[n]}$. We then also have that for every $x\in [n]$ that $\Delta(\mathcal{C}_0^x)=\Delta(\mathcal{C}_1^x)=\gamma_x=0$. And by induction we have that $G(\mathcal{C}_0^x)$ and $G(\mathcal{C}_1^x)$ have at most $2$ components each. This then clearly implies there are at most $4$ components in $G(\mathcal{C})$. If there is actually $4$ such components then both of $G(\mathcal{C}_0^x)$ and $G(\mathcal{C}_1^x)$ have $2$ components each and no component in $G(\mathcal{C})$ will cross both of $G(\mathcal{C}_0^x)$ and $G(\mathcal{C}_1^x)$. Now then this means that for every complementary pair in $\mathcal{C}$, say $A,\overline{A}\in \mathcal{C}$ that if $x\not \in A$ then $A\in \mathcal{C}_0^x$ and so lands in a family that only has elements from $\mathcal{C}_0^x$ and so if $A\in \mathcal{C}_1^x$ then it would connect the family that only had elements from $\mathcal{C}_1^x$ and the family that had $A$. This then means that $|d(\mathcal{C})|\leq |R_x|$, but as we have that both $G(\mathcal{C}_0^x)$ and $G(\mathcal{C}_1^x)$ are disconnected and $\Delta(\mathcal{C}_0^x)=\Delta(\mathcal{C}_1^x)=0$, we have that from Lemma \ref{fullshadow} that $|d(\mathcal{C}_0^x)|,|d(\mathcal{C}_1^x)|>0$ which contradicts the assumption that $\gamma_x=0$.
Now assume that $G(\mathcal{C})$ had three components, calling them $X_1,X_2,X_3$, now if we had that the number of components of $G(\mathcal{C}_0^x)$ plus $G(\mathcal{C}_1^x)$ was $3$, similar reasoning to the $4$ family case would show that $\gamma_x>0$.
Given the fact that both $G(\mathcal{C}_1^x)$ and $G(\mathcal{C}_0^x)$ have $2$ components, it must then be that one of the slices, will sit alone in a component.

Letting $M_i=\{x:\exists A_1,A_2\in {X}_i\text{ such that } x\in A_1\text{ and } x\not\in A_2\}$ and $A\in A_j$ we let $u_{i,j}=A\cap M_i$. Whenever $x\in M_i$ it is the case that the component $X_i$ will meet both $\mathcal{C}_0^x$ and $\mathcal{C}_1^x$ and by Lemma \ref{fullshadow} and the fact that $\gamma_x=0$ it is the case that there are complementary pairs between all$X_i$ and $X_j$ for $j\not=i$ and so from Lemma \ref{Variationproperties} we have that $M_i$ partition $[n]$ and $\emptyset\subsetneq u_{i,j}\subsetneq M_i$.

 We now take $S=u_{1,2}\cup u_{2,3}\cup u_{3,1}$ and show that in fact $S\not\in \Gamma(\mathcal{C})$. Now we have that for every $A\in X_1$, that $A\cap u_{2,3}=\emptyset$ and so $S\not\subseteq A$ and since all $M_i$ are disjoint, and $\emptyset\subsetneq u_{i,j}\subsetneq M_i$ we get that $A\not\subset S$. A similar argument for $X_2$ and $X_3$ then shows that $S\not\in \Gamma(\mathcal{C})$ which is a contradiction by Lemma \ref{fullshadow} and so $\mathcal{C}$ did not have $3$ connected components and so has at most $2$.
\end{proof}
From now on we assume towards a contradiction that $|G|=2^{n-1}-1$. 
We then also define $e(\mathcal{C})=\{\{A,\overline{A}\}:A,\overline{A}\in \mathcal{U_C}\setminus \mathcal{D_C}\text{ or } A,\overline{A}\in \mathcal{D_C}\setminus \mathcal{U_C}\}$ which is clearly disjoint with $\mathcal{D_C}\cap \overline{\mathcal{U_C}}$. This then implies that $|b(\mathcal{C})|\geq |P|-\frac{1}{2}|P\cap \overline{P}|+|e(\mathcal{C})|$.

From Lemma $\ref{fullshadow}$, $\ref{2comp}$ and \ref{convexblock} we get that $\Delta(\mathcal{C})=1$, $\Gamma(\mathcal{C})=2^{[n]}$,  $|s(\mathcal{C})|=|\mathcal{C}|-|G|$ and $|e(\mathcal{C})|=0$ for which we assume in all the remaining Lemmas. We also note that by Lemma \ref{convexblock}, having $m$ connected components in $\mathcal{C}=con(G)$ implies that we have at least $m$ connected components in $G$.
\begin{lemma}
It is the case that $k=3$.
\end{lemma}
\begin{proof}
  For any $x\in [n]$ assume that $\mathcal{C}_0^x$ meets $\leq 1$ family, this then means that at least $2$ of the $k$ original families would have to be represented solely by $\mathcal{C}_1^x$, yet we assumed that complementary pairs existed between every $2$ original families. A similar argument for $\mathcal{C}_1^x$ then implies a contradiction. Now given $\Delta(\mathcal{C})=1$, we have that either $\Delta(\mathcal{C}_0^x)=0$ or $\Delta(\mathcal{C}_1^x)=0$ and without loss of generality we take $\Delta(\mathcal{C}_1^x)=0$. By Lemma \ref{2comp} we have that $\mathcal{C}_1^x$ has at most $2$ connected components and so it meets at most $2$ original families, and so this forces $k=3$ again by the assumption on complementary pairs between families.
\end{proof}
\begin{lemma}
There exists some $x\in [n]$ such that $\mathcal{C}_0^x$ or $\mathcal{C}_1^x$ meets all of $\mathcal{A}_1,\mathcal{A}_2,\mathcal{A}_3$
\end{lemma}
\begin{proof}
Assume towards a contradiction that no such $x\in [n]$ exists. We have that for every such $x\in [n]$ that only $1$ of the families have sets $A_1,A_2$ such that $x\not\in A_1$ and $x\in A_2$. Let $M_i=\{x:\exists A_1,A_2\in \mathcal{A}_i\text{ such that } x\in A_1,x\not\in A_2\}$ and for $i\not=j$ and some $A\in X_j$ we let $u_{i,j}=A\cap M_i$ we then have that by Lemma \ref{Variationproperties} it is the case $M_i$ partition $[n]$, are non-empty and $\emptyset\subsetneq u_{i,j}\subsetneq M_i$. Therefore as before $S=u_{1,2}\cup u_{2,3}\cup u_{3,1}\not\in \Gamma({\mathcal{C}})$, which is a contradiction. 
\end{proof}
We then would like to chose such $x\in [n]$ where one of $\mathcal{C}_0^x$ or $\mathcal{C}_1^x$ spans $3$ components. Without loss of generality we can assume that $\mathcal{C}_0^x$ spans 3 components, as if $\mathcal{C}_1^x$ is the slice meeting all $3$ components, replacing every set in $\mathcal{A}_i$ by its complement means that $(\overline{\mathcal{C}})_0^x$ meets all 3 components, and so we relabel $\overline{\mathcal{C}}$ as $\mathcal{C}$ if necessary. We have that $\Delta(\mathcal{C}_0^x)=1$, $\Delta(\mathcal{C}_1^x)=0$ and $\gamma_x=0$.

We then have a very specific structure for our families, we have that $\mathcal{C}_0^x$ spans all three components and $\mathcal{C}_1^x$ spans $2$ of them, so relabelling if necessary we assume that $\mathcal{C}_1^x$ does not meet $\mathcal{A}_3$.

We also note the following Lemma;
\begin{lemma}
\label{nofullvariation}
If $\Delta(\mathcal{C})=1$ then it is the case that for no $y\in [n]$ can both $\mathcal{C}_0^y$ and $\mathcal{C}_1^y$ span all $3$ families
\end{lemma}
\begin{proof}
Since $\Delta(\mathcal{C})=1$, one of $\Delta(\mathcal{C}_0^y)=0$ or $\Delta(\mathcal{C}_1^y)=0$ occurs, and so without loss of generality assume that $\Delta(\mathcal{C}_0^y)=0$, then if $\Gamma(\mathcal{C}_0^y)=2^{[n]\setminus\{y\}}$ then by lemma \ref{2comp} $\mathcal{C}_0^y$ has at most $2$ components which would contradict the assumption that $\mathcal{C}_0^y$ spanned all $3$ families, and so it must be that the shadow is incomplete but then Lemma \ref{downsetint} says that $\Delta(\mathcal{C}_0^y)\geq 1$ again a contradiction.
\end{proof}

We call a family $\mathcal{A}\subseteq 2^{[n]}$ \emph{$(X,t)$-complete} if it has the property that $\{A:A\cap X=t\}\subseteq \mathcal{A} $
\begin{lemma}
\label{vary}
We take $\cS$ to be a convex family with $\cS_1,\cS_2$ being its only connected components, such that $\Delta(\cS)=0$. There exists a partition $[n]=X\sqcup Y\sqcup Z$ such that we can characterise our families as $\cS_1=\{A\subseteq [n]:A\cap X\in U_{\cS_1}$ and $A\cap Y\not\in U_{\cS_2}\}$ and $\cS_2=\{A\subseteq [n]:A\cap X\not\in U_{\cS_1}$ and $A\cap Y\in U_{\cS_2}\}$
\end{lemma}
\begin{proof}
We let $D=D_{S}$ and $U=U_{S}$ and since $\cS_i$ are convex we have that $\cS_i=D_{S_i}\cap U_{S_i}$ and since the components $\cS_i$ are disconnected we have that $D_{\cS_1}\cap U_{\cS_2}=U_{\cS_1}\cap D_{\cS_2}=\emptyset$. From Lemma \ref{fullshadow} $\Gamma(\cS)=2^{[n]}$ and so $$|D\cap \overline{U}|-|\cS|=|2^{[n]}\setminus (D\cup \overline{U})|=|(U\setminus D)\cap(\overline{D\setminus U})|=|(D\setminus U)\cap \overline{(U\setminus D)}|.$$ 
Since we have that $\Delta(\cS)=0$ we will get that $d(\cS)=|(D\setminus U)\cap \overline{(U\setminus D)}|$. We then let $K_1=D_{\cS_1}\cap \overline{U_{\cS_2}}$ and $K_2=D_{\cS_2}\cap \overline{U_{\cS_1}}$. Now $K_1\cap K_2\subseteq (D\setminus U)\cap \overline{(U\setminus D)}$ and so we get that \begin{equation*}
|d(\cS)|=|(D\setminus U)\cap \overline{(U\setminus D)}|\geq|K_1\cap K_2|\geq \frac{|K_1|| K_2|}{2^n}\geq|K_1\cap \overline{K_2}|=|d(\cS)|
\end{equation*}
This then means that $(D\setminus U)\cap \overline{(U\setminus D)}=K_1\cap K_2$. Let  $X=\{x:|(K_{1})_0^x|>|(K_1)_1^x|\}
$ and $Y=\{x:|(K_2)_0^x|>|(K_2)_1^x|\}$ and so Lemma \ref{downsetequality} says that $X\cap Y=\emptyset$, now we let $Z=[n]\setminus(X\cup Y)$. Now we clearly have that $K_1$ and $K_2$ are non empty (by $d(S)\not=\emptyset$) and $K_i\not=2^{[n]}$ as $D_{\cS_1}\cap U_{\cS_2}=U_{\cS_1}\cap D_{\cS_2}=\emptyset$. Then we can define $\mathbb{D}_1=\{A\cap X:A\in K_1\}$ and $\mathbb{D}_2=\{A\cap Y:A\in K_2\}$ which are clearly downsets. 
Now we let $\mathbb{P}_1=K_1\cap \overline{K_2}$ and $\mathbb{P}_2=K_2\cap \overline{K_1}$. Now for any $A\in \mathbb{P}_2$ it must be the case that $A\cap Y\in \mathbb{D}_2$ and $X\setminus (A\cap X)\in \mathbb{D}_1$. We therefore have that for any $A\in D_{\mathbb{P}_2}$ that $A\subseteq S\in \mathbb{P}_2$ and so $A\cap Y\subseteq S\cap Y$, and for any $S$ such that $S\cap Y\in \mathbb{D}_2$ we have that for $M=S\cup X$ we then have that $M\cap Y\in \mathbb{D}_2$ and $\emptyset \in \mathbb{D}_1$, and so we have shown that $D_{\mathbb{P}_2}=\{A:A\cap Y\in \mathbb{D}_2\}$, a similar argument shows that $U_{\mathbb{P}_1}=\{A:Y\setminus A\in \mathbb{D}_2\}$

We have that if for some $t\subseteq X$ that $t\in \mathbb{D}_1$, we then have by definition of $X$ that $\{A: A\cap X=t\}\subseteq K_1$ and so is $(X,t)$-complete, we have that such sets are present in $D_{
\cS_1}$ and $\overline{U_{\cS_2
}}$ (and so all sets with fixed intersection with $X$ being $X\setminus t$ are present in $U_{\cS_2}$). Now since $D_{\cS_1}\cap U_{\cS_2}=\emptyset$ we have that $U_{\cS_2}$ does not have sets with $A\cap X=t$ (since of course they would show up in the intersection), and a similar argument works to show that no such set with fixed intersection with $X$ being $X\setminus t$ exists in $D_{\cS_1}$ and so clearly $X\setminus t\not\in \mathbb{D}_1$.
(We note that $K_1\cap \overline{K_2}\subseteq D_1\cap U_1$ and so both $t$ and $X\setminus t$ can't be in $\mathbb{D}_1$).

Now if we have that $t\not\in \mathbb{D}_1$ and $X\setminus t\not \in \mathbb{D}_1$ then we let $A=X\setminus t$ and $\overline{A}$, which by definition does not exist in $K_1$.We then have that since $\mathbb{D}_2$ is a downset that $A\in D_{\mathbb{P}_2}$ and so $A\in D_{\cS_2}$. A similar argument shows then that
$\overline{A}\in U_{\cS_1}$. Now if $A\not\in U_{\cS_2}$ and $\overline{A}\not\in D_{\cS_1}$ then we would have that $A\in (D\setminus U)\cap \overline{(U\setminus D)}$, yet we know that $A\not\in K_1$ which contradicts the fact that $(D\setminus U)\cap \overline{(U\setminus D)}=K_1\cap K_2$. Then it is the case that either $A\in U_{\cS_2}$ or $\overline{A}\in D_{\cS_1}$. If it is the case that $\overline{A}\in D_{\cS_1}$ then we have that $D_{\cS_1}$ clearly will be $(X,t)$-complete and if $A\in U_{\cS_2}$ then $U_{\cS_2}$ will be $(X,X\setminus t)$-complete  , it is then also not possible that $A\in U_{\cS_2}$ and $\overline{A}\in D_{\cS_1}$ since then $\overline{A}\in K_1$ and so (by the previous lemma we can generate all supersets that have elements from Y and Z) $\{B:B\cap X=t\}\subseteq K_1$ which by the fact that $t\not\in \mathbb{D}_1$ this would not be possible. Now by repeating the previous argument with $A=t$ would give, then either tell us that $D_{\cS_1}$ is both $(X,t)$-complete and $(X,X\setminus t)$-complete or $U_{\cS_2}$ is both $(X,t)$-complete and $(X,X\setminus t)$-complete (since $K_1$ is neither $(X,t)$-complete nor $(X,X\setminus t)$-complete), and by the fact that $U_{\cS_2}\cap D_{\cS_1}=\emptyset$, we have that $D_{\cS_1}=2^{[n]}\setminus U_{\cS_2}$ and (by applying the same argument to $U_{\cS_1}$ and $D_{\cS_2}$) we have that $D_{\cS_2}=2^{[n]}\setminus U_{\cS_1}$. Since $\cS_1=D_{\cS_1}\cap {U_{\cS_1}}$ it is then the case that after relabelling $X$ and $Y$ the result follows.
\end{proof}

    As a short reminder we have that assuming that $|G|=2^{n-1}-1$ that for $\mathcal{C}=con(G)$ that $\Delta(\mathcal{C})=1$, $\Gamma(\mathcal{C})=2^{[n]}$,  $|s(\mathcal{C})|=|\mathcal{C}|-|G|$, $|e(\mathcal{C})|=0$ and $x\in[n]$ is such that $\mathcal{C}_0^x$ spans all $3$ families, and we let $A_3$ be the set that only sees $\mathcal{C}_0^x$. We also note that $D^x_0=D_{\mathcal{C}_1^x}\cup D_{\mathcal{C}_0^x}$, $D^x_1=D_{\mathcal{C}_1^x}$, $(\overline{U})_0^x=\overline{U_{{\mathcal{C}_1^x}}}\cup \overline{U_{\mathcal{C}_0^x}}$ and $(\overline{U})_1^x=\overline{U_{{\mathcal{C}_0^x}}}$ 

Now we note that it is the case that $\mathcal{C}_0^x\subseteq U_{\mathcal{C}_1^x}$, which can be seen as follows: if for any $A\in \mathcal{C}_0^x$ we would have that $A\in D_{\mathcal{C}_1^x}\setminus U_{\mathcal{C}_1^x}$ we have by definition that $A\cap X\not\in U_{\cS_1}$ and $A\cap Y\not\in U_{\cS_2}$, and so it is the case that $A\in D_{\cS_1}\cap D_{\cS_2}$, and since $A\in \mathcal{C}$ we would then have that $A$ would connect the original 2 families $\mathcal{A}_1$ and $\mathcal{A}_2$ (since $\cS_1$ and $\cS_2$ lived in $\cA_1$ and $\cA_2$ respectively).

 We then fix $A_1,A_2\in \mathcal{A}_3$ such that $\overline{A_1}\in \mathcal{A}_1$ and $\overline{A_2}\in \mathcal{A}_2$, which of course exist by the assumption that $G$ is non-degenerate. We then apply Lemma \ref{vary} to $\mathcal{C}_1^x$ to get a partition of $[n]\setminus \{x\}=Q$ being $X\sqcup Y \sqcup Z$ such that the $2$ connected components of $\mathcal{C}_1^x$ are labelled as $\cS_1=\{A\subseteq Q:A\cap X\in U_{\cS_1}$ and $A\cap Y\not\in U_{\cS_2}\}$ and $\cS_2=\{A\subseteq Q:A\cap X\not\in U_{\cS_1}$ and $A\cap Y\in U_{\cS_2}\}$.

We will then note some important relations.

\begin{lemma}
\label{properties}
The following are true:
\begin{enumerate}
\item[$\mathbf{1}$.] ${A_1}\cap X,X\setminus ({A_1}\cap X)\in U_{\cS_1}$ and ${A_2}\cap Y,Y\setminus ({A_2}\cap Y)\in U_{\cS_2}$ 
\item[$\mathbf{2}$.] $A_1\cap Y\in U_{\cS_2}$,$Y\setminus (A_1\cap Y)\not\in U_{\cS_2}$ and $A_2\cap X\in U_{\cS_1}$,$X\setminus (A_2\cap X)\not\in U_{\cS_1}$ 

\end{enumerate}
\end{lemma}
\begin{proof}
For $\mathbf{1}$ take $A_1\cap X=t$, we have that since for any $t\subseteq X$ we have that either $U_1$ or $D_2$ are $(X,t)$- complete, then if $A_1\cap X\not\in U_1$ then it would be the case that $D_2$ would be $(X,t)$-complete, would then imply there is some $B\in \cS_2$ such that $A_1\subseteq B$ which is not possible by the incomparability condition between families. Additionally say $D_2$ is $(X,X\setminus t)$ complete, then we would have that there is some $B\in \cS_2$ such that $\overline{A_1}\setminus\{x\}\subseteq B$ and so $\overline{A_1}\subseteq B\cup\{x\}$ which again is not possible since $\overline{A_1}$ lived in $\mathcal{A}_1$. Doing the similar operation for $A_2\cap Y$ gives the 2nd part of $\mathbf{1}$.

For $\mathbf{2}$, similar to before we take $A_1\cap Y=t$ if $t\not\in U_{\cS_2}$ then it would be the case that $D_{\cS_1}$ is $(Y,t)$-complete and so $A_1\in D_{\cS_1}$. Therefore there is some $B\in \cS_1$ where $A_1\subseteq B\subseteq B\cup\{x\}$ which contradicts the incomparability relation. Now, assume that $Y\setminus t\in U_{\cS_2}$ then it is the case that $U_{\cS_2}$ is $(Y,Y\setminus t)$-complete and so we would get that there is some $B\in \cS_2$ such that  $B\subseteq \overline{A_1}\setminus \{x\}$, and so $B\cup \{x\}\subseteq \overline{A_1}$ which again is not possible since $\overline{A_1}$ lives in $\mathcal{A}_1$. Doing the similar operation for $A_2\cap X$ gives the 2nd part of $\mathbf{2}$.
\end{proof}
\begin{lemma}
\label{existence}
We have that if $A,\overline{A}\in U_{\mathcal{C}_1^x}\setminus D_{\mathcal{C}_1^x}$ then it must be that $A$ and $\overline{A}$ live in the same connected component of $\mathcal{C}_0^x$.
\end{lemma}
\begin{proof}
We first show existence. Since we have that $\Gamma(\mathcal{C})=2^{[n]}$ we have that $A,\overline{A}\in \Gamma(\mathcal{C})$, if it was the case that $A\not\in D_{\mathcal{C}_0^x}$  then we would get that $\{A,\overline{A}\cup \{x\}\}\subseteq U\setminus D$, but we would then have that $|e(\mathcal{C})|>0$ which is a contradiction and so we have that $A\in D_{\mathcal{C}_0^x}$. Now we then have that if say $\overline{A}\not\in U_{\mathcal{C}_0^x}$ then we would have, letting $\overline{U}=I$, that ${A}\in (D_0^x\setminus D_1^x)\cap (I^x_0\setminus I_1^x)$. Now we know, by the fact that $A,\overline{A}\in U_{\mathcal{C}_1^x}\setminus D_{\mathcal{C}_1^x}$ that neither $A\cup \{x\}\in C$ nor $\overline{A}\cup \{x\}\in {C}$ and so we would get that $\gamma_x>0$ which is a contradiction. The same argument for $\overline{A}$ shows that $\overline{A}\in \mathcal{C}_0^x$ as well.

Now if we had that $A,\overline{A}$ lived in different connected components of $\mathcal{C}_0^x$ then we would get that again $\gamma_x>0$ as neither of $A\cup \{x\},\overline{A}\cup \{x\}$ are in $\mathcal{C}$ by the fact that $A,\overline{A}\in U_{\mathcal{C}_1^x}\setminus D_{\mathcal{C}_1^x}$.
\end{proof}
\begin{lemma}
\label{finalcontradiction}
For any $A,\overline{A}\in U_{\mathcal{C}_1^x}\setminus D_{\mathcal{C}_1^x}$ it is the case that $A,\overline{A}\not\in con(\mathcal{A}_3)$
\end{lemma}
\begin{proof}
Assume that $A\in con(\mathcal{A}_3)$, then by Lemma \ref{existence} we have that $A,\overline{A}\in con(\mathcal{A}_3)$ and so $\mathcal{A}_3$ varies on $Q$. Now since it is the case that $A,\overline{A}\in U_{\mathcal{C}_1^x}\setminus D_{\mathcal{C}_1^x}$, we have that $A,\overline{A}\in U_{\mathcal{S}_1}\cap U_{\mathcal{S}_2}$ and so $\mathcal{S}_1$ varies on all of $X$ and $\mathcal{S}_2$ varies on all of $Y$. Since $\Delta(\mathcal{C})=1$ by Lemma \ref{nofullvariation} we can not have that all families vary on any given element of $[n]$ and since $Z$ from Lemma \ref{vary} is defined as free we have that both $\mathcal{S}_1$ and $\mathcal{S}_2$ vary on $Z$ but since $con\mathcal{A}_3$ varies on the full $Q$, it is the case that $Z=\emptyset$. 
Additionally Lemma \ref{nofullvariation} implies that $con(\mathcal{A}_1)$ is constant on $Y$ and $con(\mathcal{A}_2)$ is constant on $X$. 
By Lemma $\ref{vary}$ there exists sets in $\mathcal{S}_1$ that are disjoint from $Y$ and sets in $\mathcal{S}_2$ that are disjoint with $X$ and so we have that for every $T\in con(\mathcal{A}_1)$ that $T\cap Y=\emptyset$ and every $T\in con(\mathcal{A}_2)$ that $T\cap X=\emptyset$ (by the fact that the intersections must be constant). 
Now taking $A_1\in \mathcal{A}_3$ such that $\overline{A_1}\in \mathcal{A}_1$, we then get that $Y\subseteq A_1$, but since $\mathcal{C}_0^x$ is present in $\mathcal{A}_2$, one can take $K\in \mathcal{C}_0^x\cap con(\mathcal{A}_2)$ and as $K\cap X=\emptyset $ and it is the case that $Z=\emptyset$ and $x\not\in K$ that $K\subseteq Y$, but we also know that $Y\subseteq A_1$ contradicting the incomparability condition between families.
\end{proof}

\begin{proof}[Proof of Theorem \ref{3rd}]
Take $A=(A_1\cap X)\cup (A_2\cap Y)$ and $\overline{A}$ with respect to $Q$. 
 Now we have that by Lemma \ref{properties} that $A_1\cap X\in U_{\cS_1}$, $A_2\cap Y\in U_{\cS_2}$ and $\overline{A_1}\cap X\in U_{\cS_1}$, $\overline{A_2}\cap Y\in U_{\cS_2}$. And so we have that $A,\overline{A}\in U_{\mathcal{C}_1^x}\setminus D_{\mathcal{C}_1^x}$. By Lemma \ref{existence} we have that $A,\overline{A}\in \mathcal{C}_0^x$. Without loss of generality assume that $A,\overline{A}\in con(\mathcal{A}_1)$, we then have that $con(\mathcal{A}_1)$ varies on all of $[n]\setminus \{x\}$, and so also on all of $Y$, and by Lemma \ref{properties} the family $con(\mathcal{A}_2)$ also varies on all of $Y$. Since $\Delta(\mathcal{C})=1$ by Lemma \ref{nofullvariation} it can not be the case that all 3 convexified families vary on the same element of $[n]$, and so all sets in $con(\mathcal{A}_3)$ must have a fixed intersection with $Y$, namely $A_2\cap Y=A_1\cap Y$. Yet from Lemma \ref{properties} we have that $Y\setminus{A}_1\not\in U_{\mathcal{S}_2}$ and $Y\setminus A_2\in U_{\mathcal{S}_2}$ which is a contradiction. The same argument shows why $A,\overline{A}$ can not live in $con(\mathcal{A}_2)$ and so they must both live in $con(\mathcal{A}_3)$. Since if $A,\overline{A}\in \mathcal{A}_3$ implies that $A,\overline{A}\in con(\mathcal{A}_3)$, by Lemma \ref{finalcontradiction} this is a contradiction, and so the result follows.
\end{proof}
\section{Conclusion}
We will conclude the paper by providing an example of $3$ families that are non-degenerate for which their sizes sum to $2^{n-1}-4$.

Fix $n\ge6$, put $Q=[n]\setminus[4]$ and $q=2^{|Q|}=2^{n-4}$.

For families $\mathcal R$ and $\mathcal T$ let
$\mathcal R\vee\mathcal T=\{R\cup T:R\in\mathcal R,\,T\in\mathcal T\}$.
\[
\begin{aligned}
 \mathcal A_1&=\{3,4,34\}\vee\bigl(2^Q\setminus\{\emptyset\}\bigr)\ \cup\
        \{13,14,134\}\vee\bigl(2^Q\setminus\{\emptyset,Q\}\bigr)\ \cup\ \{134\},\\
 \mathcal A_2&=\{23,24,234,123,124\},\\
 \mathcal A_3&=\{2,12\}\vee\bigl(2^Q\setminus\{\emptyset\}\bigr)\ \cup\ \bigl\{\{1\}\cup Q\bigr\}.
\end{aligned}
\]
So $|\mathcal A_1| = 3(2^{n-4}-1) +3 (2^{n-4}-2) + 1 = 3\cdot 2^{n-3} -8$, $|\mathcal A_2| = 5$ and $|\mathcal A_3| = 2(2^{n-4}-1) + 1 = 2^{n-3} - 1$. The families are uncomplemented and pairwise incomparable with
$|\mathcal A_1 \cup \mathcal A_2 \cup \mathcal A_3|=2^{n-1}-4$.
\printbibliography
\end{document}